\documentclass[12pt]{amsart}  
\usepackage{amssymb,amsmath,amsthm,amscd}
\usepackage[all]{xy}

\usepackage[hypertex]{hyperref}

\numberwithin{equation}{section}

\newtheoremstyle{fancy1}{10pt}{10pt}{\itshape}{12pt}{\textsc\bgroup}{.\egroup}{8pt}{
}
\newtheoremstyle{fancy2}{10pt}{10pt}{}{12pt}{\itshape}{.}{8pt}{ }

\theoremstyle{fancy1}

\newtheorem{cor}[equation]{Corollary}
\newtheorem{lem}[equation]{Lemma}
\newtheorem{prop}[equation]{Proposition}
\newtheorem{thm}[equation]{Theorem}
\newtheorem{problem}{Problem}

\newtheorem*{main*}{Theorem}

\newtheorem*{cor*}{Corollary}

\newtheorem*{problem*}{Problem}

\theoremstyle{fancy2}

\newtheorem*{rem*}{Remark}

\newcommand{\cref}[1]{Corollary~\ref{#1}}

\newcommand{\lref}[1]{Lemma~\ref{#1}}
\newcommand{\pref}[1]{Proposition~\ref{#1}}

\newcommand{\tref}[1]{Theorem~\ref{#1}}

\newcommand{\gt}{\theta}

\newcommand{\gs}{\sigma}

\newcommand{\e}{\epsilon}

\newcommand{\RP}{\mathbb{R\mkern1mu P}}
\newcommand{\CP}{\mathbb{C\mkern1mu P}}
\newcommand{\HP}{\mathbb{H\mkern1mu P}}
\newcommand{\CaP}{\mathrm{Ca}\mathbb{\mkern1mu P}^2}
\newcommand{\Sph}{\mathbb{S}}

\newcommand{\C}{{\mathbb{C}}}
\newcommand{\R}{{\mathbb{R}}}
\newcommand{\Z}{{\mathbb{Z}}}
 
\newcommand{\N}{{\mathbb{N}}}
\newcommand{\QH}{{\mathbb{H}}}

\newcommand{\F}{\ensuremath{\operatorname{F}}}
\newcommand{\G}{\ensuremath{\operatorname{G}}}

\newcommand{\SO}{\ensuremath{\operatorname{SO}}}

\renewcommand{\O}{\ensuremath{\operatorname{O}}}
\newcommand{\Sp}{\ensuremath{\operatorname{Sp}}}
\newcommand{\U}{\ensuremath{\operatorname{U}}}
\newcommand{\SU}{\ensuremath{\operatorname{SU}}}
\newcommand{\Spin}{\ensuremath{\operatorname{Spin}}}

\newcommand{\T}{\ensuremath{\operatorname{T}}}
\renewcommand{\S}{\ensuremath{\operatorname{S}}}

\newcommand{\K}{\ensuremath{\operatorname{K}}}
\renewcommand{\L}{\ensuremath{\operatorname{L}}}

\newcommand{\fg}{{\mathfrak{g}}}
\newcommand{\fk}{{\mathfrak{k}}}
\newcommand{\fh}{{\mathfrak{h}}}
\newcommand{\fm}{{\mathfrak{m}}}

\newcommand{\fp}{{\mathfrak{p}}}

\newcommand{\fsu}{{\mathfrak{su}}}
\newcommand{\fu}{{\mathfrak{u}}}
\newcommand{\ft}{{\mathfrak{t}}}

\def\con#1=#2(#3){#1 \equiv #2 \bmod{#3}}

\newcommand{\tr}{\ensuremath{\operatorname{tr}}}
\newcommand{\diag}{\ensuremath{\operatorname{diag}}}

\newcommand{\rank}{\ensuremath{\operatorname{rk}}}

\renewcommand{\Im}{\ensuremath{\operatorname{Im}}}
\newcommand{\Ad}{\ensuremath{\operatorname{Ad}}}
\newcommand{\ad}{\ensuremath{\operatorname{ad}}}

\renewcommand{\sec}{\ensuremath{\operatorname{sec}}}

\DeclareMathOperator{\cohom}{cohom}

\DeclareMathOperator{\Id}{Id}

\newcommand{\no}{\noindent}
\newcommand{\co}{{cohomogeneity}}
\newcommand{\coo}{{cohomogeneity one}}
\newcommand{\com}{{cohomogeneity one manifold}}
\newcommand{\coms}{{cohomogeneity one manifolds}}
\newcommand{\coa}{{cohomogeneity one action}}
\newcommand{\coas}{{cohomogeneity one actions}}

\renewcommand{\sc}{sectional curvature}
\newcommand{\cu}{curvature}
\newcommand{\nn}{non-negative}
\newcommand{\nnc}{non-negative curvature}
\newcommand{\nnsc}{non-negative sectional curvature}
\newcommand{\pc}{positive curvature}
\newcommand{\psc}{positive sectional curvature}
\newcommand{\pcu}{positively curved}
\newcommand{\qp}{quasi positive}
\newcommand{\ap}{almost positive}

\begin{document}

\title{Examples of Riemannian Manifolds with non-negative sectional curvature}

\author{Wolfgang Ziller}
\address{University of Pennsylvania\\
   Philadelphia, PA 19104}
\email{wziller@math.upenn.edu}

\maketitle

Manifolds with non-negative sectional curvature have been of
interest since the beginning of global Riemannian geometry, as
illustrated by the theorems of Bonnet-Myers, Synge, and the sphere
theorem. Some of the oldest conjectures in global Riemannian
geometry, as for example the  Hopf conjecture on $\Sph^2\times
\Sph^2$, also fit into this subject.

For non-negatively curved manifolds, there are  a number of
obstruction theorems known,  see Section 1 below and  the survey by
Burkhard Wilking in this volume. It is somewhat surprising that the
only further obstructions to positive curvature are given by the
classical Bonnet-Myers and Synge theorems on the fundamental group.

Although there are many examples with \nnc, they all come from two
basic constructions, apart from taking products. One is taking an
isometric quotient of a compact Lie group equipped with a
biinvariant metric and another a gluing procedure due to Cheeger and
recently significantly generalized by Grove-Ziller. The latter
examples include a rich class of manifolds, and give rise to \nnc\
on many exotic 7-spheres. On the other hand, known manifolds with
positive sectional curvature are very rare, and are all given by
quotients of compact Lie groups, and, apart from the classical rank
one symmetric spaces, only exist in dimension below 25.

Due to this lack of knowledge, it is therefore of importance to
discuss and understand known examples and find new ones. In this
survey we will concentrate on the description of known examples,
although the last section also contains suggestions where to look
for new ones. The techniques used to construct them are fairly
simple. In addition to the above, the main tool is a deformation
described by Cheeger that, when applied to non-negatively curved
manifolds, tends to increase curvature. Such Cheeger deformations
can be considered as the unifying theme of this survey.
 We can thus be fairly explicit in
the proof of the existence of all known examples which should make
the basic material  understandable at an advanced graduate student
level. It is  the hope of this author that it will thus encourage
others to study this beautiful subject. This survey originated in
the  Rudolph Lipschitz lecture series the author gave at the
University of Bonn in 2001 and various courses taught at the
University of Pennsylvania.

\section{General structure theorems}
\label{general}

To put the examples discussed in this survey into context, we
mention the main structure theorems and conjectures in this subject.
 See the survey by
Burkhard Wilking in this volume for further information.

\begin{itemize}
\item
(Gromov) If $M^n$ is a compact manifold with $\sec\ge 0$, then there
exists a universal constant $c(n)$ such that $b_i(M^n,F)\le c(n)$
for all $i$ and any field of coefficients $F$. Furthermore, the
fundamental group has a generating set with at most $c(n)$ elements.

\item (Cheeger-Gromoll)
If $M^n$ is a compact manifold that admits a metric with
non-negative sectional curvature,  then there exists an abelian
subgroup of $\pi_1(M^n)$ with finite index.

\item (Lichnerowicz-Hitchin) The obstructions to positive scalar curvature
imply that a compact spin manifold with $\hat{A}(M)\ne 0$ or
$\alpha(M)\ne 0$ does not admit a metric with non-negative sectional
curvature. This holds in particular for the unique exotic sphere in
dimension 9.

\item
(Cheeger-Gromoll) If $M^n$ is a non-compact manifold with a complete
metric with $\sec\ge 0$, then there exists a totally geodesic
compact submanifold $S^k$, called the soul, such that $M^n$ is
diffeomorphic to the normal bundle of $S^k$.
\end{itemize}

Surprisingly, for positive curvature one  has in addition only the
classical  obstructions:

\begin{itemize}
\item
(Bonnet-Myers) A manifold which admits a metric with positive
curvature has finite fundamental group.
\item
(Synge) An even dimensional manifold with positive curvature has
fundamental group $0$, if orientable, and $\Z_2$, if non-orientable.
In odd dimensions a positively curved manifold is orientable.
\end{itemize}

If we allow ourselves to add an upper as well as a lower bound on
the sectional curvature it is convenient to introduce what is called
the {\it pinching constant} which is defined as $\delta = \min \sec
/ \max\sec$. One then has the following recognition and finiteness
theorems:

\begin{itemize}
\item
(Berger-Klingenberg) If $M^n$ is a compact simply connected manifold
with $\delta\ge \frac 1 4 $, then $M$ is either homeomorphic to
$\Sph^n$ or isometric to $\CP^n$, $\HP^n$ or $\CaP$ with their
standard Fubini metric.
\item
(Cheeger) Given a positive constant $\e$, there are only finitely
many diffeomorphism types of compact simply connected manifolds
$M^{2n}$ with $\delta\ge \e$.
\item (Fang-Rong,Petrunin-Tuschmann) Given a positive constant $\e$, there are only finitely many
diffeomorphism types of compact   manifolds $M^{2n+1}$ with
$\pi_1(M)=\pi_2(M)=0$ and $\delta \ge \e$.
\end{itemize}

We finally mention some conjectures.

\begin{itemize}
\item
(Hopf) There exists no metric with positive sectional curvature on
$\Sph^2\times\Sph^2$. More generally, there are no positively curved
metrics on the product of two compact manifolds, or on a symmetric
space of rank at least two.
\item (Hopf) A compact manifold with $\sec\ge 0$ has non-negative Euler
characteristic. An even dimensional manifold with positive curvature
has positive Euler characteristic.
\item (Bott) A
compact simply connected manifold $M$ with $\sec\ge 0$ is {\it
elliptic}, i.e.,  the sequence of Betti numbers of the loop space of
$M$
 grows at
 most polynomially for every field of coefficients.
\end{itemize}
The latter conjecture, and its many consequences, were discussed in
the literature for the first time in \cite[\!\!'82]{GH}.

\section{Compact examples with non-negative curvature}

There are two natural constructions to produce new metrics with
non-negative curvature from given ones. If $M_1$ and $M_2$ are
endowed with metrics of \nnc, the product metric on $M_1\times M_2$
clearly has \nnc. The second construction is by taking quotients, or
more generally by considering Riemannian submersions.

Recall that  if  $M$ and $B$ are  two Riemannian manifolds, then a
smooth map $\pi\colon M\to B$ is called a {\it Riemannian
submersion} if $\pi_*$ is an isometry on horizontal vectors, i.e.,
on vectors orthogonal to the fibers. For such submersions one has
the O'Neill formula:
$$
\sec_B(\pi_* x,\pi_* y)=\sec_M(x,y)+\tfrac 3 4 ||\; [X,Y]^v\; ||^2 ,
$$
where $x,y$ are orthonormal horizontal vectors, i.e., orthogonal to
the fibers, $X,Y$ are horizontal vector fields extending $x,y$, and
$[X,Y]^v$ denotes the vertical part of $[X,Y]$, i.e., the component
tangent to the fiber. If $M$ has non-negative curvature, so does
$B$, and one can hope that in some cases $B$ is even positively
curved. The most basic examples of Riemannian submersions are given
by taking quotients $\pi\colon M\to M/G$ where $G$ is a compact Lie
group acting freely and isometrically on $M$.  We often call the
induced metric on $M/G$ the `quotient' metric.

Before we describe a third method, let us first  recall some
standard ways of putting metrics on homogeneous spaces. If a compact
Lie group $G$ acts transitively on $M$ and $p\in M$, we can write
$M=G/H$ where $H$ is the isotropy group at $p$. We will often fix a
biinvariant metric $Q$ on $\fg$, the Lie algebra of $G$. Note that
$\ad X \colon Y\to [X,Y]$ is then a skew symmetric endomorphism with
respect to $Q$. Thus the $Q$-orthogonal decomposition $\fg=\fh+\fm$
satisfies $[\fh,\fh]\subset\fh$ and $[\fh,\fm]\subset\fm$. The more
restrictive condition $[\fm,\fm]\subset\fh$ corresponds to the case
where the metric $Q$ induces a locally symmetric metric on $G/H$. We
identify $\fm$ with $T_pM$ via action fields: $X\in\fm \to X^*(p)$
where $X^*(q)=\frac{d}{dt}_{|t=0}\exp(tX)q$. The action of $H$ on
$T_pM$ is then identified with the action of $\Ad H$ on $\fm$. A
metric on $G/H$, invariant under the action of $G$, corresponds to
an inner product on $\fm\cong T_pM$ which is $Ad(H)$-invariant. This
inner product can be expressed as $Q(PX,Y)$ where $X,Y\in\fm$ and
$P\colon\fm\to\fm$ is a $Q$-symmetric endomorphism.

A third method that produces new \nn ly curved metrics from a given
one is obtained via a  {\it Cheeger deformation}. This process was
first used by M.Berger who considered metrics on spheres, shrunk in
the direction of the Hopf fibration, to produce odd dimensional
manifolds with small injectivity radius and positively pinched
curvature. A systematic general description was given in
\cite[\!\!'73]{Che}. Let $(M,g)$ be a Riemannian manifold and $G$ a
Lie group acting by isometries on $M$. We then consider the
Riemannian submersion
\begin{equation*}
\pi\colon M\times G \to M \quad \colon\quad (p,g)\to g^{-1}p\ .
\end{equation*}
This can also be viewed as a quotient construction via the action
\begin{equation*}\label{Cheeger1}
\tilde{g}\star(p,g)=(\tilde{g}p,\tilde{g}g) \text{\ or simply \ }
M=(M\times G)/\triangle G .
\end{equation*}
We can thus start with a \nn ly curved metric $g$ on $M$, take a
product with a biinvariant metric, and then the quotient metric
defines a new metric on $M$. To describe this process as a
deformation, fix a biinvariant metric $Q$ on $\fg$ and let $g_t$ be
the metric obtained as a quotient of the product metric $g+\frac 1 t
Q$ on $M\times G$.  Tangent to the orbit $Gp=G/G_p$, we write the
original metric  as above in the form $Q(PX,Y)$, where
$X,Y\in\fm_p$, with $\fm_p$ the orthogonal complement of the Lie
algebra of $G_p$. The symmetric endomorphism $P\colon\fm_p\to\fm_p$
is changed into a new symmetric endomorphism $P_t$ describing $g_t$
in terms of $Q$ and we claim:
\begin{equation}\label{Cheeger2}
P_t=(P^{-1}+t\Id)^{-1} .
\end{equation}
To see this, observe that $\pi_*(X^*,Y)=X^*-Y^*$. Thus a horizontal
lift of $X\in\fm_p\cong T_p(G/G_p)\subset T_pM$, under the
Riemannian submersion $\pi$, is equal to
\begin{equation*}\label{lift}
(P^{-1}(P^{-1}+t\Id)^{-1}X^*(p),-t(P^{-1}+t\Id)^{-1}X)\in\T_pM\times\fg
,
\end{equation*}
and the length squared of this vector is
\begin{align*}
Q((P^{-1}+t\Id)^{-1}X,P^{-1}(P^{-1}+t\Id)^{-1}X)&+\frac 1 t
Q(t(P^{-1}+t\Id)^{-1}X,t(P^{-1}+t\Id)^{-1}X)\\
&=Q((P^{-1}+t\Id)^{-1}X,X) .
\end{align*}
Orthogonal  to the orbit $Gp$, the metric is  unchanged since a
horizontal lift of $X\in\fm_p^\perp\subset T_pM$ is equal to
$(X^*(p),0)$.

This process can be considered as a deformation, since we obtain the
original metric $g$ when $t=0$. If $\lambda_i$ are the eigenvalues
of  $P$, the eigenvalues of $P_t$ are
$\frac{\lambda_i}{1+t\lambda_i}$, i.e. the metric is shrunk in the
direction of the orbits. This process will in general destroy
symmetries, although the group $G$ clearly still acts by isometries
induced by right multiplication in $M\times G$. We can thus also
consider iterated Cheeger deformations corresponding to a chain of
subgroups $H_1\subset \dots \subset H_k\subset G$.

The deformation $g\to g_t$ tends to improve curvature. If the
original metric $g$ has \nnc, the deformed metric does also by
O'Neill's formula. If $X,Y$ span a 0-curvature 2-plane of $g$, then,
by considering the $G$-components of vectors in $M\times G$, its
curvature becomes positive in the metric $g_t$ if
$[PX_\fm,PY_\fm]\ne 0$, where we have associated to $X\in T_pM$ a
vector  $X_{\fm}\in\fm_p$ such that $X_{\fm}^*(p)$ is the component
of $X$ in the orbit direction. Although this will not be needed in
this survey, one finds a detailed study of the basic properties of
this deformation in M\"uter \cite[\!\!'87]{Mu}. We mention here only
that, if we let $C_t=P^{-1}P_t$ on $\fm_p$ and
${C_t}_{\vert\fm_p^\perp}=\Id$  be the symmetric endomorphism that
expresses $g_t$ in terms of $g=g_0$, then
$\sec_{g_t}(C_t^{-1}X,C_t^{-1}Y)>0$ for $t>0$ unless
$\sec_g(X,Y)=0$, $[PX_{\fm},PY_{\fm}]=0$ and $d\omega_Z(X,Y)=0$ for
all $Z\in\fg$, where $\omega_Z$ is the one form dual to the Killing
vector field $Z$.  Thus the 0-curvature planes tend to ``move" with
$C_t^{-1}$. Furthermore, 2-planes which are tangent to a totally
geodesic flat 2-torus, and which contain a vector orthogonal to the
$G$ orbit, remain flat.

\bigskip

As a starting point for finding examples, one considers compact Lie
groups $G$ endowed with a biinvariant metric since their curvature
satisfies:
$$
\sec_G(x,y)=\frac 1 4 ||\, [x,y]\, ||^2\ge 0 \text{ for } x,y\in\fg
\text{ orthonormal} .
$$

Combining this fact  with O'Neill's formula, we obtain non-negative
curvature on every isometric quotient of a compact Lie group. In
particular, all homogeneous spaces $G/H$, where $H$ is a closed
subgroup of $G$, have metrics with non-negative curvature. Since the
identity component of the isometry group of a compact Lie group $G$,
endowed with a biinvariant metric, consists of left and right
translations, it is natural to generalize the class of homogeneous
manifolds to what are called {\it biquotients}. Consider $H\subset
G\times G$ and define an action of $H$ on $G$ by
$$
h\star g = h_1 g h_2^{-1} \ , \quad \text{ where } \quad
h=(h_1,h_2)\in H.
$$
The action is free if and only if $h_1$ conjugate to $h_2$, for
$(h_1,h_2)\in H$, implies that $h_1=h_2=e$. If this is the case, the
quotient is a manifold, which we denote by $G/\!/H$ and is called a
biquotient. If $H=L\times K\subset G\times G$, we will also write
$L\backslash G/K$.  Thus we obtain:
\begin{thm}
A biinvariant metric on $G$ induces a metric with non-negative
sectional curvature on every homogeneous space $G/H$ and every
biquotient $G/\!/H$.
\end{thm}


The first time where biquotients were considered in geometry, was in
\cite[\!\!'75]{GM}, where it was shown that an exotic 7-sphere
admits \nnc. To describe this example, consider the action of the
unit quaternions $\Sp(1)$ on the symplectic group $ \Sp(2)$ via:
$$
q\star A = \diag(q,q)A\diag(q,1)^{-1} \qquad q\in\Sp(1)\ ,\
A\in\Sp(2).
$$
This action is clearly free and we have:
\begin{thm}[Gromoll-Meyer]\label{GM}
The \nn ly curved \ manifold $\Sp(2)/\!/\Sp(1)$ is homeomorphic, but
not diffeomorphic, to $\Sph^7$.
\end{thm}

 In order to prove this, one observes that
$\Sp(2)/\!/\Sp(1)\Sp(1)=\diag(q,q)A\diag(r,1)^{-1}, $ $ \ q,r\in
\Sp(1)$, is diffeomorphic to $\Sph^4$ since the action of
$\diag(q,q)$ on $\Sp(2)/\diag(r,1)=\Sph^7$ is the Hopf action by
$\S^3$. Thus $\Sp(2)/\!/\Sp(1)$ can be considered as an $\Sph^3$
bundle over $\Sph^4$. One then identifies which sphere bundle it
represents by using Milnor's description of certain exotic 7-spheres
as $\Sph^3$ bundles over $\Sph^4$. We also point out that  in
\cite[\!\!\!'02]{To}, and independently in \cite[\!\!'04]{KZ}, it
was shown that the only exotic sphere which can be written as a
biquotient is the Gromoll-Meyer sphere.

\smallskip

Another special class of \nn ly curved examples were constructed in
\cite[\!\!'73]{Che}:
\begin{thm}[Cheeger]\label{rankone}
The connected sum of any two rank one symmetric spaces carries a
metric with \nnsc.
\end{thm}
In \cite[\!\!\!'02]{To} it was shown that some of these  Cheeger
manifolds, but not all, can be viewed as biquotients as well.

\bigskip

The gluing construction used in order to prove \tref{rankone}, was
recently significantly generalized to what are called {\it \coo\
manifolds}. Recall that if $G$ is a Lie group that acts on a
manifold $M$, the cohomogeneity of the action is defined as $\cohom
(M,G)=\dim M/G$. Thus an action with $\cohom(M,G)=0$ is an action
that is transitive, i.e., the manifold is a homogeneous spaces.
Cohomogeneity one manifolds can thus be considered as the next
simplest kind of group actions to study. They are also special among
all group actions since, as we will see, the manifold can be
reconstructed from its isotropy groups. The geometry and topology of
homogeneous spaces is fairly well understood by now, whereas this is
not yet the case for \coo\ manifolds. One should point out though,
that this class of manifolds does not contain the homogeneous spaces
as a subset. In fact only very few homogeneous spaces carry a \coo\
action.

\smallskip

Let $G$ be a compact group acting by \coo\ on a compact manifold
$M$. Since $M/G$ is one dimensional, it is either a circle $\S^1$,
or an interval $I$. In the first case all $G$ orbits are principal
and $\pi\colon M\to M/G=\S^1$ is a fiber bundle with fiber a
principal orbit $G/H$, and thus the fundamental group  is infinite.
One also easily sees that such fiber bundles carry a $G$ invariant
metric with \nnc. In the second more interesting case there are
precisely two nonprincipal $G$-orbits corresponding to the endpoints
of $I$, and $M$ is decomposed as the union of two tubular
neighborhoods of these nonprincipal orbits, with common boundary a
principal orbit. Let $\ell_-$ and $\ell_+$ be the codimension of the
nonprincipal orbits.  We have the following existence theorem
(\cite[\!\!'00]{GZ1}):
\begin{thm}[Grove-Ziller]\label{codim2}
A compact \coo\ $G$-manifold with $\ell_\pm \le 2$ has a
$G$-invariant metric with \nnsc.
\end{thm}
One easily sees that $\ell_\pm>1$ if $M$ is simply connected.
Although the assumption on the codimensions seems rather special, it
turns out that the class of \coms\ with $\ell_\pm=2$ is surprisingly
rich. An immediate application is:
\begin{cor}
Each of the 4 oriented diffeomorphism types of $\RP^5$'s carries a
metric with \nnsc.
\end{cor}
This follows since $\Sph^5$ carries (non-linear) \coa s by
$\SO(2)\SO(3)$, discovered by G.Calabi, with codimension 2 singular
orbits. They are a special case of the action on Kervaire spheres
described below. The involution in $\SO(2)$ acts freely and, using
surgery theory, one shows  that any  one of the exotic $\RP^5$'s can
be obtained in this fashion, see \cite[\!\!'61]{Lo}.

\smallskip

In \cite[\!\!'00]{GZ1} it was also conjectured that \tref{codim2}
holds without any assumption on the codimensions. This turns out to
be false. One has (\cite[\!\!\! '06]{GVWZ}):
\begin{thm}[Grove-Verdiani-Wilking-Ziller]\label{nonex}
For each pair $(\ell_-,\ell_+)$ with $(\ell_-,\ell_+)\ne (2,2)$ and
$\ell_\pm \ge 2$ there exist infinitely many \coo\ $G$-manifolds
that do not carry a $G$ invariant metric with \nnsc.
\end{thm}

The most interesting example in this Theorem are the Kervaire
spheres, which are the only  exotic spheres that can carry a \coo\
action (\cite[\!\!'03]{St}). They are described as  a $2n-1$
dimensional Brieskorn varietiy:
$$ z_0^d + z_1^2 +\cdots z_n^2=0 \quad , \quad |z_0|^2 + \cdots |z_n|^2 =1 .$$
It carries
 a \co\ one action by \SO(2)\SO(n)
defined by $(e^{i\gt},A)(z_0,\cdots ,
z_n)=(e^{2i\gt}z_0,e^{id\gt}A(z_1,\cdots,z_n)^t\, )\, $ whose
codimensions are $(\ell_-,\ell_+)=(2,n)$. For $n$ odd and $d$ odd,
they are homeomorphic to spheres, and are exotic spheres if $\con
2n-1=1(8)$. If $n\ge 4, d\ge 3$ one shows that there exists no
$G$-invariant metric with \nnc.

To prove \tref{nonex}, one needs to describe the set of all
$G$-invariant metrics explicitly. They depend on a finite collection
of functions, 6 in the case of the Kervaire spheres, which need to
satisfy certain smoothness conditions at the endpoint. For each
choice of 2-planes tangent to the principal orbit, one obtains
differential inequalities on these functions and their first
derivatives from the Gauss equations. By a suitable choice of
2-planes one obtains upper and lower bounds on the first derivatives
 which contradict the smoothness conditions at one of the singular
orbits.

We mention that in the case of Ricci curvature one has the positive
result that every \com\ carries an invariant metric with
non-negative Ricci curvature, and with positive Ricci curvature if
and only if the fundamental group is finite (\cite[\!\!'02]{GZ2}).

\smallskip

To discuss the proof of \tref{codim2} and some of its applications,
we first recall the basic structure of \coas. We will only consider
the most interesting case, where $M/G = I$ and let $\pi\colon M\to
M/G$ be the projection. In order to  make the description more
explicit, we choose an arbitrary but fixed $G$-invariant Riemannian
metric on $M$, normalized so that with the induced metric, $M/G =
[-1,1] $. Fix a point $x_0 \in \pi^{-1}(0)$ and let $c : [-1,1] \to
M$ be a geodesic orthogonal to the orbit through $x_0$, and hence to
all orbits, and parameterized such that $\pi\circ c = id_{[-1,1]}$.
 Let
$B_\pm = \pi^{-1}(\pm 1) = G\cdot x_\pm$ be the two nonprincipal
orbits, where $x_\pm = c(\pm 1)$. It then follows that
 $c : [2n-1,2n+1] \to M$, $n\in \Z$ are
minimal geodesics between the two nonprincipal orbits $B_\pm$ since
$G$ acts transitively on the set of all geodesics orthogonal to the
orbits. Let $K_\pm = G_{x_\pm}$ be the isotropy groups at $x_\pm$
and $H = G_{x_0} = G_{c(t)}$, $ -1<t<1 $, the principal isotropy
group. By the slice theorem, we have the following description of
the tubular neighborhoods $D(B_-)=\pi^{-1}([-1,0])$ and
$D(B_+)=\pi^{-1}([0,1])$ of the nonprincipal orbits $B_\pm =
G/K_\pm$ :
\begin{equation*}
\label{discbundle} D(B_{\pm}) = G\times_{K_{\pm}}D^{\ell_{\pm}} ,
\end{equation*}
\no where $D^{\ell_{\pm}}$ is the normal (unit) disk to $B_\pm$ at
$x_\pm$. Here the action of $K_{\pm} $ on $G\times D^{\ell_{\pm}}$
is given by $k\star (g,p)=(gk^{-1},kp)$ where $k$ acts on
$D^{\ell_{\pm}}$ via the slice representation.  Hence we have the
decomposition
\begin{equation*}
\label{decomp} M=D(B_-) \textstyle{\: \cup}_E \: D(B_+) \; ,
\end{equation*}
\no where $E=\pi^{-1}(0) = G\cdot x_0 = G/H$ is a principal orbit
which is canonically identified with the boundaries $\partial
D(B_\pm) = G\times_{K_{\pm}}\Sph^{\ell_{\pm}-1}$, via the maps $G\to
G\times \Sph^{\ell_{\pm}-1}$ , $g\to (g,\mp\dot{c}(\pm1))$. Note
also that $\partial D^{\ell_{\pm}} = \Sph^{\ell_{\pm}-1} = K_\pm /H$
since the boundary of the tubular neighborhoods must be a $G$ orbit
and hence $\partial D^{\ell_{\pm}}$ is a $K_\pm$ orbit.  All in all
we see that we can recover $M$ from $G$ and the subgroups $H$ and
$K_\pm$.  We caution though that the isotropy types, i.e., the
conjugacy classes of the isotropy groups $K_\pm $ and $ H$ do not
determine $M$.

\smallskip

An important fact about \coas\ is that there is a converse to the
above construction. Suppose $G$ is a compact Lie group and $H\subset
K_\pm \subset G$ are closed subgroups, which we sometimes denote by
$H\subset \{K_-,K_+\}\subset G$. Assume furthermore  that $K_\pm/H =
\Sph^{\ell_{\pm}-1}$ are spheres. It is well known  that a
transitive action of a compact Lie group $K$ on a sphere $\Sph^\ell$
is conjugate to a linear action and is determined by its isotropy
group $H\subset K$. We can thus assume that $K_\pm$ acts linearly on
$\Sph^{\ell_\pm}$ with isotropy group $H$ at $p_\pm\in
\Sph^{\ell_\pm-1}$ and define a manifold
\begin{equation*}
\label{manifold} M = G\times_{K_{-}}D^{\ell_{-}} \cup_{G/H}
G\times_{K_{+}}D^{\ell_{+}} ,
\end{equation*}
\no where we glue the two boundaries by sending $[g,p_-]$ to
$[g,p_+]$.
 $G$ acts on $M$   via $g^*[g,p]=[g^*g,p]$ on each half and one easily
checks that it  has isotropy groups $K_\pm $ at $[e,0]$ and $ H$ at
$[e,p_0]$ and is thus \coo.
\bigskip

\tref{codim2} clearly follows from the following geometric result by
gluing two such metrics on the tubular neighborhoods $D(B_\pm)$
along their common boundary $G/H$.

\begin{prop}[]\label{metricdisk}
Let $H\subset K\subset G$ be Lie groups with $K/H=\Sph^{1}=\partial
D^{2} $ and fix a biinvariant metric $Q$ on $G$.
  On the disc bundle
 $G\times_{K} D^{2}$  there exists
a $G$-invariant metric  with \nnsc, which is a product near the
boundary $G\times_K\Sph^1=G/H$ with  metric on $G/H$ induced by $Q$.
\end{prop}

The crucial ingredient in the proof of \pref{metricdisk} is the
following result about left invariant metrics.

\begin{lem}\label{leftinv}
Let $G$ be a compact Lie group and  $\fk\subset\fg$ an abelian
subalgebra. Consider the left invariant metric on $G$ whose value at
$T_eG=\fg$ is given by
$$Q_t=tQ_{\mid\fk}+ Q_{\mid \fk^\perp} \ ,
$$
where $Q$ is a biinvariant metric on $G$. Then $Q_t$ has \nnsc\ as
long as $t\le 4/3$.
\end{lem}
\begin{proof}
The curvature formula for a left invariant metric $\langle X ,Y
\rangle=Q(PX,Y)$ is given by (see e.g. \cite[\!\!'99]{Pu}):
\begin{align*}\label{curv}
\langle R(X,Y)Y,X)\rangle =& \tfrac 1 2
Q([PX,Y]+[X,PY],[X,Y])-\tfrac 3 4 Q(P[X,Y],[X,Y]) \\\notag
 &+ Q(B(X,Y),P^{-1}B(X,Y)) - Q(B(X,X),P^{-1}B(Y,Y))\ ,
\end{align*}
where $B (X,Y) = \tfrac{1}{2} ([X,PY] - [PX,Y])$.

In our case, let $X=A+R,Y=B+S$ with $A,B\in\fm=\fk^\perp$ and
$R,S\in\fk$ and hence  $P(A+R)=A+tR$. We can now split up the
expressions into  components in direction of $\fm$ and of $\fk$. A
computation shows that the $\fm$ component is given by
$$\tfrac{1}{4}\left\|[A,B]_\fm + t[X,B] + t [A,Y] \right\|^2_Q\ge
0 \ ,$$ where we have used the bi-invariance of $Q$ and the Jacobi
identity to show that \newline $\langle [X,B],[A,Y]\rangle=\langle
[X,A],[Y,B]\rangle$. On the other hand,  the $\fk$ component is
given by
$$
\|[A,B]_{\fk}\|_Q^2-\tfrac 3 4 t\|[A,B]_{\fk}\|_Q^2=(1-\tfrac 3 4
t)\|[A,B]_{\fk}\|_Q^2 \ ,
$$
which is non-negative as long as $t\le \frac 4 3 $.
\end{proof}

{\it Proof of \pref{metricdisk} } We have inclusions $H\subset
K\subset G$ with $K/H=\Sph^1$ and define $Q$-orthogonal
decompositions $\fg=\fk+\fm$ and $\fk=\fh+\fp$. As usual, we
identify the tangent spaces $ T_{(H)}K/H\cong \fp$ and
$T_{(H)}G/H\cong \fp+\fm$. Since $\fp$ is  one dimensional,
\lref{leftinv} implies that the left invariant metric on $G$ defined
by $Q_a=aQ_{|\fp}+Q_{|\fh+\fm}$ has \nnc\ as long as $a\le 4/3$.
Since $[\fp,\fp]=0$  and $[\fh,\fp]\subset\fp$, the subalgebra $\fp$
is an ideal of $\fk$ and hence $\Ad K$ invariant, and thus $Q_a$ is
right $K$-invariant as well. In addition we choose a metric
$g_f=dt^2+f(t)^2d\theta^2$ on $D^2$ which is clearly invariant under
the action of $K$ on $D^2$ and has \nnc\ if $f$ is concave. The
product metric $Q_a+g_f$ on $G\times D^2$ thus induces a \nn ly
curved metric $g_{a,f}$ on the homogeneous disk bundle $G\times_K
D^2$. We now claim that given $1<a\le 4/3$, we can choose $f$ such
that $g_{a,f}$ is a product near the boundary with metric on $G/H$
induced by $Q$. To see this, consider the Riemannian submersion
$G\times (K/H) \to G\times_K K/H\cong G/H$ where we endow
$K/H=\Sph^1$ with the metric of a circle of radius $f(t)$. The
induced metric on $G/H$ is the metric $g_{a,f}$ restricted to  the
boundary of a tube of radius $t$. We compute this metric as in the
case of a Cheeger deformation \eqref{Cheeger2}. If $2\pi s_0$ is the
length of the circle $K/H$ in the metric $Q_{|\fp}$, the metric on
$K/H$ is given by $(f/s_0)^2Q$ and  it follows that the metric on
$G/H$ is given by $Q$ on $\fm$ and by
$\frac{a}{1+a(f/s_0)^{-2}}Q=\frac{f^2a}{f^2+as_0^2}Q$ on $\fp$.
Hence we obtain the desired metric by choosing a concave function
$f$ and a $t_0$ such that $f^2(t)=\frac{as_0^2}{a-1}$, for $t\ge
t_0$.
  \qed

\begin{rem*} We can view this construction as a ``scaling up, scaling
down" procedure. The natural metric on $G\times_{K}D$ induced by a
biinvariant metric $Q$ on $G$ shrinks the metric on the boundary
$G/H$ in the direction of $K/H$, as in the case of a Cheeger
deformation. This needs to be compensated by scaling the metric $Q$
up in the direction of $\fp$ in order to recover the metric $Q$ on
$G/H$. This explains the difficulty of proving \pref{metricdisk} for
$\ell>2$ since left invariant metrics as in \lref{leftinv} in
general have some sectional curvature positive when $t>1$. In fact
we believe:

\begin{problem}
Let  $G$ be a compact simple Lie group and $K$ a non-abelian
subgroup. Show that a metric $Q_t$ as in \lref{leftinv} has some
negative sectional curvatures for any $t>1$.
\end{problem}
\end{rem*}

Nevertheless, it is  possible that there are other special
homogeneous disk bundles for which \pref{metricdisk} holds with
codimension $\ell>2$.

\smallskip

As was observed by B.Wilking, \tref{codim2} can be generalized to
the situation where the homogeneous orbits are replaced by
biquotients. In other words, if $K_\pm\subset G\times G$ acts freely
on $G$ and $H\subset K_\pm$ with $K_\pm/H=\Sph^1$, then the
resulting manifold carries a metric with \nnc. This follows by
applying \tref{codim2} to the \com\ $H\subset \{K_-,K_+\}\subset
G\times G$ and then dividing by $\Delta G\subset G\times G$ on the
left, which acts freely by assumption.

\smallskip

We now apply this result to some concrete \coms\ in order to prove:
\begin{thm}[Grove-Ziller]\label{principals4}
Every principal $\SO(k)$ bundle $P$ over $\Sph^4$ carries a \coa\ by
$\SO(3)\times\SO(k)$ with codimension two principal orbits and hence
an invariant metric with non-negative curvature.
\end{thm}

Thus, by O'Neill's formula, every associated bundle
$P\times_{\SO(k)}X$ with $X$ a non-negatively curved manifold on
which $\SO(k)$ acts by isometries, also carries a non-negatively
curved metric.
\begin{cor}
Every sphere bundle over $\Sph^4$ carries a metric with \nnsc.
\end{cor}

Of particular interest are  $\Sph^3$ bundles over $\Sph^4$ since
Milnor discovered the first exotic spheres among these manifolds. It
implies:
\begin{cor}
Of the $14$ (unoriented) exotic 7-spheres, $10$ carry a metric with
non-negative curvature.
\end{cor}
The group of exotic spheres, under the group operation of connected
sums, is isomorphic to $\Z_{28}$, but a change of orientation
corresponds to taking an inverse. It is not known whether  the
remaining 4 exotic spheres  carry \nnc\ metrics as well.

\bigskip

{\it Proof of \tref{principals4} }:  Let the \com\ $P_{r,s}$ be
given by the isotropy groups:
$$
H=\triangle Q\subset \{ (e^{ir\gt},e^{i\gt})\cdot H \, ,
(e^{js\gt},e^{j\gt}\cdot H\}\subset \S^3\times \S^3,
$$
\no where $\triangle Q=\{\pm(1,1),\pm(i,i),\pm(j,j),\pm(k,k),\}$ is
the quaternion group and $ e^{ir\gt}=\cos(r\gt)+i\sin(r\gt)$ is an
embedding of a circle into the unit quaternions $\S^3$. In order for
$H$ to be a subgroup of $K_\pm$, we need to assume that $ \con
r,s=1(4)$. We then have   $K_\pm/H=\S^1$ and thus \tref{codim2}
implies that $P_{r,s}$ carries an $\S^3\times \S^3$ invariant metric
with \nnc. The subgroup $\S^3=\S^3\times \{e\}\subset
\S^3\times\S^3$ acts freely on $P_{r,s}$ since its isotropy groups
are the intersection of $\S^3\times \{e\}$ with $K_\pm$ and $H$,
which by construction are trivial. We now claim that $P_{r,s}/\S^3$
is $\Sph^4$. To see this, observe that the second $\S^3$ factor
induces a \coa\ with group diagram $Q\subset \{e^{i\gt}\cdot
Q,e^{j\gt}\cdot Q\}\subset \S^3$ on the quotient. The element
$-1\in\S^3$ acts trivially and the effective version of the action
has isotropy groups
$\Z_2\times\Z_2\subset\{\S(\O(2)\O(1)),\S(\O(1)\O(2))\}\subset
\SO(3)$.  But there is a well known  linear action by      $\SO(3)$
on $\Sph^4$ given by conjugation on the set of $3\times 3$ symmetric
real matrices with trace $0$. Since every matrix is conjugate to a
diagonal one, it follows that the two singular orbits are given by
symmetric matrices with two equal eigenvalues, positive for one and
negative for the other, and the principal orbits by matrices with 3
distinct eigenvalues. One now easily checks that the isotropy groups
are the same as for the above action and hence $P_{r,s}/\S^3$  is
equivariantly diffeomorphic to $\Sph^4$.

Thus $P_{r,s}$ can be viewed as an $\S^3$ principal bundle over
$\Sph^4$. These are classified by an integer $k$, namely the Euler
class of the bundle evaluated on a fixed orientation class of
$\Sph^4$. To recognize which bundle it is, one observes that the
Gysin sequence implies  $H^4(E,\Z)=\Z_{|k|}$ for such a bundle. For
a \com\ one can compute the cohomology groups by using
Meyer-Vietoris on the decomposition into the disk bundles
$D(B_\pm)$. The disk bundles are homotopy equivalent to $G/K_\pm$
and their intersection to $G/H$. Using well known methods for
computing the cohomology groups of homogeneous spaces one shows that
$H^4(P_{r,s},\Z)$ is a cyclic group of order $(r^2-s^2)/8$. But for
$ \con r,s =1(4)$ the values of $(r^2-s^2)/8$ can take on any
integer. Thus every $\S^3$ principle bundle over $\Sph^4$ is of the
form $P_{r,s}$ for some $r,s$. Since every $\SO(3)$ principal bundle
over $\Sph^4$ is spin, i.e., has a lift to an $\S^3$ principal
bundle, this implies \tref{principals4} for $k=3$. The case of $k=4$
one obtains by repeating the above argument for
$G=\S^3\times\S^3\times\S^3$ with $\K_\pm$ again one dimensional and
identity component of say $K_-$ equal to
$(e^{ir_1\gt},e^{ir_2\gt},e^{i\gt})$ with $r_i\in\Z$. For principal
bundles $P$ over $\Sph^4$ with $k>4$ it is well known that their
structure  reduces to a  $\SO(4)$. Thus there exists an
$\SO(4)$-principle bundle $P'$ with $P=P'\times_{\SO(4)}\SO(k)$ on
which $\SO(k)$ acts on the right. Hence the lift of $\SO(3)$ to $P'$
also lifts to $P$ and commutes with $\SO(k)$. \qed

\bigskip

We finally indicate how the proof of Cheeger's \tref{rankone} fits
into the above framework. Of the connected sums considered in his
theorem, only $\CP^n \# - \CP^n$ admits a \coa.  But a similar idea
as in the proof of \pref{metricdisk} applies to all cases. A rank
one projective space $M^n$  with a small ball removed, is
diffeomorphic to the disk bundle of the canonical line bundle over
the projective space of one dimension lower. This bundle is  a
homogeneous disk bundle with boundary diffeomorphic to a sphere. One
now uses the same ``scaling up, scaling down" method as in the proof
of \pref{metricdisk} to show that these disk bundles have a metric
with non-negative curvature  which is a product near the boundary
and has constant curvature one on the boundary. One can then glue
together any two rank one projective spaces along this boundary.

\bigskip

The  methods described  in the proof  of \tref{principals4} can also
be applied to other 4-manifolds as base (\cite[\!\!'07]{GZ3}):
\begin{thm}[Grove-Ziller]\label{principalcp2}
Every principal $\SO(k)$ bundle $P$ over $\CP^2$ which is not spin,
i.e., $w_2(P)\ne 0$, carries a \coa\ with codimension two principal
orbits and hence an invariant metric with non-negative curvature.
Thus, so does every associated sphere bundle.
\end{thm}
To prove this, one uses the linear \coa\ on $\CP^2$ given by
$\SO(3)\subset\SU(3)$, which one easily verifies has  group diagram
$\Z_2\subset\{\S(\O(2)\O(1)),\SO(2)\}\subset\SO(3)$ and constructs a
group diagram with $G=\S^3\times\S^3$ as above, but with $H=\{(\pm
1,\pm1), $ $( \pm i,\pm i)\}$. The topological considerations needed
to identify what bundle the \com\ represents, are significantly more
difficult. It also raises  the following general question, which the
above examples show is important in the context of \coms.

\smallskip

{\it Given a principal $\L$-bundle $P \to M$ over a $\G$-manifold
$M$. When does  the action of $\G$ on $M$ have a {\it commuting
lift}, i.e., a lift to an action of $\G$, or possibly a cover of
$\G$, on the total space $P$, such that the lift commutes with $L$.
}

\smallskip

This problem has been studied extensively. However, apart from the
general result that every action of a semi simple group admits a
commuting lift to the total space of every principal circle or more
generally torus bundle \cite[\!\!'61]{PS}, the results seem to be
difficult to apply in concrete cases. For a \com\ $M_G$ with
isotropy groups $H\subset \{K_-,K_+\}\subset G$, one has a natural
description of the lifts to an $L$-principal bundle over $M_G$ in
terms of the isotropy groups. Simply choose embeddings of $K_\pm$
into $L\times G$ such that they agree on $H$ and are given in  the
second component by the original embeddings into $G$. The action by
$L\times \{e\}$ is then clearly free, and the quotient is $M_G$
since the  induced $G$ action has the same isotropy groups. As long
as one allows the action of $G$ on $M_G$ to be ineffective, all
lifts are described in this fashion. The difficulty is then to
decide what the isomorphism type of this $L$-principle bundle is.

\tref{principals4} and \tref{principalcp2} can be restated as saying
that the linear actions of $\SO(3)$ on $\Sph^4$ and $\CP^2$ have a
commuting lift to every principal $\SO(k)$ bundle, respectively
principal $\SO(k)$ bundle which is not spin. In \cite[\!\!'07]{GZ3}
one finds a classification of which \coas\ on simply connected
4-manifolds $M^4$ have a commuting lift to a given principal
$\SO(k)$ bundle over $M^4$. In particular, it is shown that in the
spin case the action of $\SO(3)$ on $\CP^2$ only lifts to half of
all $\SO(3)$ principal bundles. This shows the limitations of our
principle bundle method which finds metrics on their total space
with $\sec\ge 0$..

\vspace{5pt}

A particularly interesting case of the above Problem are $\SO(k)$
principal bundles over $\Sph^k$ since cohomogeneity one actions  on
spheres are numerous and have been classified in \cite[\!\!'71]{HL}.
\begin{problem}
Which \coas\ on $\Sph^n$ admit a commuting lift to a given $\SO(k)$
principal bundle over $\Sph^n$?
\end{problem}
An answer to this question could potentially produce further sphere
bundles over spheres, and hence possibly higher dimensional exotic
spheres, which carry metrics with \nnc.

\smallskip

In light of the existence \tref{codim2} and the non-existence
\tref{nonex}, it is natural to pose the following somewhat vague but
important:

\begin{problem}
How large is the class of \com s  that  admit an invariant metric
with \nnc?
\end{problem}

 Are there other obstructions, and how
strong are the obstructions developed in the proof of \tref{nonex}?
As far as existence is concerned, one would need to understand how
to put non-negative curvature on \coms\ without making the middle
totally geodesic. In \cite[\!\!'02]{Sc} Schwachh\"ofer showed that
for the adjoint action of $\SU(3)$ on $\Sph^7\subset \fsu(3)$ there
exist no invariant metric with \nnc\ such that the middle is totally
geodesic (for any homogeneous metric on the principal orbit!). But
there of course exists an invariant metric with \pc.
\smallskip

We end this section with the following natural problem. Many
examples are obtained by taking  a quotient of a compact Lie group,
equipped with a left invariant metric with $\sec\ge 0$,  by a group
acting by isometries. It thus seems to be important to know what all
such metrics look like.

\begin{problem}
Classify all left invariant metrics with non-negative sectional
curvature on compact Lie groups.
\end{problem}

Surprisingly, the only examples known so far are obtained by
combining the following: Cheeger deformations of a biinvariant
metric along a subgroup $K$, i.e. the metric on $G=G\times_K K$
induced by $Q+\frac 1 t Q$.  If the subgroup is 3-dimensional, we
can more generally consider  the  metric on $G$ induced by $Q+\frac
1 t g$ where $g$ is a left invariant metric on $K$ with positive
curvature. Finally, we can scale a biinvariant metric up in the
direction of an abelian subalgebra as in \lref{leftinv}.
 The only Lie groups
where a complete answer is known, are  $\SU(2)$ and $\U(2)$, see
\cite[\!\!'06]{BT}, and with partial results for $\SO(4)$,
\cite[\!\!'06]{TH}. In the latter paper it was also shown, as
another application of  Cheeger deformations, that every left
invariant \nn ly curved metric $g$ on a compact Lie group $G$ can be
connected by an ``inverse linear" path of \nn ly curved left
invariant metrics to a fixed biinvariant metric $Q$. Indeed, in the
Cheeger deformation \eqref{Cheeger2} applied to the right action of
$G$ on itself, we can let $t\to\infty$ and then the rescaled metric
$tg_t\to Q$ since the eigenvalues of $tg_t$ in terms of $Q$ are
$\frac{t\lambda_i}{1+t\lambda_i}$. Thus the main interest lies in
deciding what derivatives are allowed at $Q$ for an inverse linear
path of  left invariant metrics with $\sec\ge 0$. This approach is
discussed in detail in \cite[\!\!'06]{TH}.

\bigskip

\centerline{\it Topology of \nn ly curved manifolds}

\bigskip

For the following we assume that our manifolds are compact and
simply connected. Recall the Bott conjecture which states that a \nn
ly curved manifold is elliptic. Even rationally elliptic, i.e.,
where the  condition on the Betti numbers of the loop space is only
assumed for rational coefficients, already has strong consequences.
By Sullivan's theory of minimal models in rational homotopy theory,
rationally elliptic is equivalent to saying that there are only
finitely many homotopy groups which are not finite.
 Rationally elliptic implies that  the sum of the Betti number of $M^n$ is at
most $2^n$, which is the optimal upper bound in Gromov's Betti
number theorem. Furthermore, the Euler characteristic is
non-negative (one half of the Hopf conjecture), and positive if and
only if the odd Betti numbers are 0. Thus it is natural to
conjecture that an even dimensional manifold with \pc\ has vanishing
odd Betti numbers. See \cite[\!\!'82]{GH} where geometric
consequences for rationally elliptic as well as for the remaining
class of simply connected so-called rationally hyperbolic manifolds
were first discussed and treated.

In dimension four, rationally elliptic manifolds are homeomorphic to
one of the known examples with \nnc, i.e., one of $\Sph^4$, $\CP^2$,
$\Sph^2\times\Sph^2$ or $\CP^2 \#\pm \CP^2$. It is natural to
conjecture that they are indeed diffeomorphic, and that only the
first two can admit \pc. In \cite{PP} it was shown that an elliptic
5-manifold is diffeomorphic to one of the known examples with \nnc,
i.e., one of $\Sph^5$, $\SU(3)/\SO(3)$, $\Sph^3\times\Sph^2$ or the
non-trivial $\Sph^3$ bundle over $\Sph^2$. Thus the Bott conjecture
in dimension 5 states that a \nn ly curved manifold is diffeomorphic
to one of these models, and it is natural to conjecture that only
the first admits \pc.

We remark that in dimension two a \nn ly curved manifold is
diffeomorphic to $\Sph^2$ by Gauss-Bonnet and in dimension three to
$\Sph^3$ by Hamilton's theorem \cite[\!\!'82]{Ha}.

\smallskip

We now describe some topological properties of the known examples
with \nnc. A homogeneous manifold $M$ is 2-connected iff $M=G/H$
with $G$ and $H$ semisimple and hence there are only finitely many
such manifolds in each dimension. If it is not 2-connected, $M$ is a
torus bundle over a 2-connected one. In \cite[\!\!'02]{To} it was
shown that both statements also hold for biquotients.

On the other hand, the class of biquotient manifolds is
significantly larger than the class of homogeneous spaces. For
example, in \cite[\!'03]{To2} it was shown that there exist
infinitely many 6-dimensional biquotients of the form
$(\,\S^3)^3/\!/(\,\S^1)^3$ with non-isomorphic rational cohomology
rings. On the other hand, compact simply connected homogeneous
spaces in dimension 6 are either diffeomorphic to a product of rank
one symmetric spaces or to the Wallach manifold $\SU(3)/T^2$.

\smallskip

The class of \com s, including associated bundles and quotients, is
again much larger than both. Indeed, there are infinitely many
homotopy types of 2-connected \com s, starting in dimension 7, since
all $\Sph^3$ bundles over $\Sph^4$ admit \nnc. In Section 6 one also
finds an infinite family of 7 dimensional \com s,  depending on 4
arbitrary integers, which are 2-connected and have singular orbits
of codimension two.

\smallskip

In \cite[\!\!'04]{DT} it was shown that there are infinitely many
non-negatively curved manifolds lying in distinct cobordism classes.
One starts with one of the principal $\SO(3)$ bundles $P$ over
$\Sph^4$ in \tref{principals4} and considers the associated bundle
$P\times_{\SO(3)}\CP^2$, where $\SO(3)\subset\SU(3)$ acts linearly
on $\CP^2$. It clearly has \nnc\ and a computation of the Pontryagin
classes shows that they have distinct Pontryagin numbers and hence
lie in different cobordism groups. On the other hand, this is not
possible for homogeneous spaces and biquotients since circle bundles
are the boundary of the associated disk bundle and hence have
vanishing Pontryagin numbers.

\smallskip

In \cite[\!\!'07]{Ho} C.Hoelscher classified  compact simply
connected \com s of dimension at most seven. In dimension 4 this was
done in \cite[\!\!'86]{Pa} (dimension 2 and 3 being trivial). But in
dimension 5, 6 and 7 there are many  \co\ actions with singular
orbits of codimension 2. In dimension  7 there are also some new
families whose codimensions are not both two, where it is not known
if they carry invariant metrics with \nnc. On the other hand, one
also has the exotic Kervaire spheres in dimension 7, which by
\tref{nonex} does not admit an invariant metric with \nnc.
\bigskip

\section{Non-Compact examples with non-negative curvature}

For non-compact manifolds one has the well known Soul Theorem
(\cite[\!\!'72]{CG}):

\begin{thm}[Cheeger-Gromoll]
If $M^n$ is a non-compact manifold with a  complete metric with
$\sec\ge 0$, then there is a totally geodesic compact submanifold
$S^k$ such that $M^n$ is diffeomorphic to the normal bundle of
$S^k$.
\end{thm}
The submanifold $S^k$ is called the soul of $M^n$. A major open
problem in this part of the subject is hence:

\begin{problem}
What vector bundles over compact manifolds with non-negative
curvature admit a complete metric with \nnsc?
\end{problem}

This is particularly interesting for vector bundles over spheres.
Any homogeneous vector bundle $G\times_KV$, where $K$ acts
orthogonally on a vector space $V$, clearly admits such a metric by
O'Neill's formula. Thus $T\Sph^n=\SO(n+1)\times_{\SO(n)}\R^n$ also
does. Every vector bundle over $\Sph^n$, $n=1,2,3$ is a homogeneous
vector bundle and hence carries \nnc. In \cite[\!\!'78]{Ri} Rigas
showed that every vector bundle of $\Sph^n$ is stably, i.e., after
taking the connected sum with a  trivial bundle of sufficiently
large dimension, a homogeneous vector bundle and hence carries
non-negative curvature.

As a consequence of \tref{principals4} and \tref{principalcp2} one
obtains non-negative curvature on the vector bundles
$P\times_{\SO(k)}\R^k$ associated to the principal bundles $P$:
\begin{cor}[Grove-Ziller]
Every vector bundle over $\Sph^4$, and every vector bundle over
$\CP^2$ which is not spin, carries a  complete metric with \nnsc.
\end{cor}

This class of vector bundles is quite large since they are
classified by one arbitrary integer when the fiber dimension is
three and by 2 if the fiber dimension is four.

\smallskip

As far as vector bundles are concerned over the remaining known
4-manifolds with non-negative curvature, i.e. $\Sph^2\times\Sph^2$
and $\CP^2 \# \pm \CP^2$, most of them also admit \nnc\ since their
structure group reduces to a torus and circle bundles over these
manifolds are known to admit \nnc, see
\cite[\!\!'95]{Ya},\cite[\!\!'02]{To}. For vector bundles over
$\Sph^n, n>4$, one  knows that all vector bundles over $\Sph^5$, and
most of the vector bundles over $\Sph^7$, admit \nnc\
(\cite[\!\!'00]{GZ1}). But in both cases there are only finitely
many such bundles.

\vspace{5pt}

If the base does not have finite fundamental group, there are
obstructions to the existence of complete metrics with \nnc\ due to
{\"O}zaydin-Walschap \cite[\!\!'94]{OW}, in the case where the soul
is flat, and Belegradek-Kapovitch \cite[\!\!'01]{BK1} ,
\cite[\!\!'03]{BK2} in general. The simplest examples are:

\begin{thm}
 Every orientable vector bundle
over $\T^2$ or $\Sph^3\times\Sph^1$  with non-negative curvature is
trivial.

\end{thm}

In \cite[\!\!'01]{BK1} , \cite[\!\!'03]{BK2} the authors give many
more examples of vector bundles over $C\times T^k$ with C compact
and simply connected and $k\ge 1$ which do not admit \nnc. For
example, if $k\ge 4$, there exist infinitely many vector bundles
over $C\times T^k$ of every rank $\ge 2$, whose total space do not
admit any complete metric with \nnc. No obstructions are known when
the base is simply connected.

\smallskip

Although it is known that for a given metric on $M^n$ any two souls
are isometric, $M$ can have two distinct non-negatively curved
metrics with souls that are not even homeomorphic. In fact
Belegradek \cite[\!\!'03]{Bel} proved:

\begin{thm}[Belegradek]
For each $n\ge 5$, there exist infinitely many complete Riemannian
metrics on $\Sph^3\times\Sph^4\times\R^n$ with $\sec\ge 0$ and
pairwise non-homeomorphic souls.
\end{thm}

To prove this, consider the principal $\SO(3)$ bundle $P_k\to
\Sph^4$ corresponding to $k\in \pi_3(\SO(3))\cong\Z$ and let
$E_k^n=P_k\times_{\SO(3)}\R^n$ and
$S_k^n=P_k\times_{\SO(3)}\Sph^{n-1}$ be the associated vector bundle
and sphere bundle coming from the standard inclusion
$\SO(3)\subset\SO(n)$. Then the bundle $\Delta^*(S_k^4\times
E_{-k}^n)$, where $\Delta\colon \Sph^4\to \Sph^4\times \Sph^4$ is
the diagonal embedding, can be regarded as a bundle over $\Sph^4$
associated to the principal $\SO(3)\times \SO(3)$ bundle
$\Delta^*(P_k^4\times P_{-k}^n)$, which by \tref{principals4}
carries an invariant metric with \nnc. On the other hand, it can
also be regarded as an n-dimensional vector bundle over $S_k^4$ and
its soul, since it is  an associated vector bundle, is equal to the
0-section $S_k^4$. Now one uses surgery theory to show that for
$\con  k= k' (12)$ and $n\ge 5$, the manifold $\Delta^*(S_k^4\times
E_{-k}^n)$ is diffeomorphic to $\Sph^3\times\Sph^4\times\R^n$ and a
computation of the Pontryagin classes shows that $S_k^4$ is
homeomorphic to $S_{k'}^4$ if and only if $k = \pm k'$.

\bigskip

See \cite[\!\!'03]{Bel} and \cite[\!\!'05]{KPT1} for further
examples of this type.

\section{Known Examples with positive curvature}

Known examples with positive curvature are surprisingly rare. What
is even more surprising is that they are all obtained as quotients
of a compact Lie group equipped with a biinvariant or a Cheeger
deformation of a biinvariant metric divided by a group of
isometries. One may view the following theorem as an explanation of
why it is so difficult to find new examples (\cite[\!\!'06]{Wi3}).

\begin{thm}[Wilking]\label{cohom} If $M^n$ admits a positively curved metric
with an  isometric action of cohomogeneity $k\ge 1$ with
$n>18(k+1)^2$, then $M$ is homotopy equivalent to a rank one
symmetric  space.
\end{thm}

Thus,  for any new examples, the larger the dimension, the bigger
the cohomogeneity. This may increase the difficulty of computing its
curvature tensor and estimating the sectional curvature. In fact,
known examples exist only in low dimensions. They consist of certain
homogeneous spaces in dimensions $6,7,12,13$ and $24$ due to Berger
\cite[\!\!'61]{Be1}, Wallach \cite[\!\!'72]{Wa}, and Aloff-Wallach
\cite[\!\!'75]{AW}, and of biquotients in dimensions $6,7$ and $13$
due to Eschenburg \cite[\!\!'82]{E1},\cite[\!\!'84]{E2} and Bazaikin
\cite[\!\!'96]{Ba1}. The purpose of this section is to discuss these
examples.

\smallskip

The main ingredient for all known examples  is the following Cheeger
deformation  of a fixed biinvariant metric $Q$ on $G$, of a type we
already considered in \lref{leftinv}. Let $K\subset G$ be a closed
Lie subgroup with Lie algebras $\fk\subset\fg$ and $\fg=\fk+\fm$ a
$Q$-orthogonal decomposition. Recall that $(G,K)$ is a symmetric
pair if $K$ is, up to components, the fixed point set of an
involutive automorphism. For our purposes, the property that
$[\fm,\fm]\subset\fk$ is all that is needed, and is equivalent to
being a symmetric pair if $G/K$ is simply connected.

 For the biinvariant metric $Q$ a
0-curvature 2-plane is characterized by $[X,Y]=0$. The following
deformation thus decreases the set of 0-curvature 2-planes
(\cite[\!\!'84]{E2}).
\begin{lem}[Eschenburg]\label{sym}
Let $Q_t$ be a left invariant metric on $G$ defined by $Q_t=tQ_{\mid
\fk} + Q_{\mid \fm}$. Then $\sec_{Q_t}\ge 0$ as long as  $t\le 1$.
If we assume in addition  that $(G,K)$ is a symmetric pair, $X,Y$
span a 0-curvature 2-plane of $g_t$, for $t<1$, if and only if
$[X,Y]=[X_{\fk},Y_{\fk}]=[X_{\fm},Y_{\fm}]=0$.
\end{lem}
\begin{proof}
The metric $Q_t$ can be viewed as a Cheeger deformation as in
\eqref{Cheeger2} with respect to the right action of $K$ on $G$ and
hence has non-negative curvature for $t\le 1$. As we saw, the metric
$Q+\frac 1 s Q$ on $G\times K$ induces a metric of the form $Q_t$
with $t=\frac{1}{s+1}<1$ and the horizontal lift of
$X=X_{\fk}+X_{\fm}\in\fk+\fm=\fg$ is equal to
$\bar{X}=(X_{\fm}+\frac{1}{1+s}X_{\fk} , -\frac{s}{1+s}X_{\fk})\in
\fg+\fk$. Since the metric on $G\times K$ is biinvariant, a
horizontal 2-plane spanned by $\bar{X},\bar{Y}$ has 0 curvature if
and only if $[\bar{X},\bar{Y}]=0$. Since the O'Neill tensor is also
given in terms of Lie brackets, the same is true for the 2-plane
spanned by $X,Y\in\fg$. If $G/K$ is a symmetric pair, we have
$[\fm,\fm]\subset \fk$, which,  together with $[\fk,\fk]\subset\fk $
and $[\fk,\fm]\subset\fm$,   easily implies the claim.
\end{proof}

Given Lie subgroups $H\subset K \subset G$,  we define  a
homogeneous fibration
\begin{equation*}\label{homog}
K/H \longrightarrow G/H \longrightarrow G/K\quad\colon\quad gH\to gK
.
\end{equation*}
Using  the $Q$-orthogonal decompositions $\fg=\fk+\fm$ and
$\fk=\fh+\fp$, we can identify the tangent spaces $\fp\cong
T_{(H)}K/H\ , \ \fm\cong T_{(K)}G/K$ and $\fp+\fm\cong T_{(H})G/H$.
 In terms of these identifications, we define a one parameter family of homogeneous metrics on $G/H$ by
\begin{equation*}\label{homogmet}
g_t=tQ_{\mid \fp} + Q_{\mid \fm},
\end{equation*}
which scales the fibers of the homogeneous fibrations by $t$. Notice
that they can also be viewed as a Cheeger deformation of the metric
$Q$ on $G/H$ in direction of the left action of $K$ on $G/H$. It is
natural to ask, if one has such a metric with \pc\ on the base and
on the fiber, when does $g_t$ have \pc. A partial answer to this
question is given by (\cite[\!\!'72]{Wa}):
\begin{prop}[Wallach]\label{fib}
Given a homogeneous fibration as above, assume that:
\begin{itemize}
\item[(a)] The base $(G,K)$ is a compact symmetric pair of rank one.
\item[(b)] The metric on the fiber $K/H$ induced by $Q$ has
positive curvature.
\item[(c)] For any non-zero vectors $X\in \fp$ and $Y\in\fm$ we have
$[X,Y]\ne 0$.
\end{itemize}
Then the metric $g_t$ with $t<1$ has positive sectional curvature.
\end{prop}
\begin{proof} The interpretation as a Cheeger deformation implies
that $\sec_{g_t}\ge 0$ for $t\le 1$. If we define the left invariant
metric $Q_t$ on $G$ by $Q_t=tQ_{\mid \fk} + Q_{\mid \fm}$, the
projection $G\to G/H$ is a Riemannian submersion with respect to the
metrics $Q_t$ and $g_t$. Thus, if $X,Y\in\fp +\fm\cong T_H G/H$ span
a 0-curvature 2-plane of $g_t$,  they span a 0-curvature 2-plane of
$Q_t$ as well and hence
$[X,Y]=[X_{\fk},Y_{\fk}]=[X_{\fm},Y_{\fm}]=0$ by \lref{sym}. The
vectors $X_{\fk},Y_{\fk}\in\fp$ can be viewed as spanning a 2-plane
of the fiber $K/H$ and since it is assumed to have positive
curvature, $X_{\fk},Y_{\fk}$ must be linearly dependent. Similarly,
since the base has \pc, $X_{\fm},Y_{\fm}$ are linearly dependent.
Hence we can find a new basis of this plane with $X\in\fp$ and
$Y\in\fm$. But now condition (c) implies that $[X,Y]=0$ is
impossible and thus $\sec_{g_t}>0$.
\end{proof}

The condition that $[X,Y]\ne 0$ is equivalent to the positivity of
the  curvature of the 2-plane  spanned by $X,Y$, i.e. the {\it
vertizontal} sectional curvatures. This condition is the fatness
condition we will discuss in Section 6.

\bigskip

\centerline{\it Homogeneous Examples with \pc.}

\bigskip

Homogeneous spaces which admit a homogeneous metric with \pc\ have
been classified by Wallach in even dimensions (\cite[\!\!'72]{Wa})
and by B\'erard-Bergery in odd dimensions (\cite[\!\!'76]{BB2}). We
now describe these examples, leaving out the compact rank one
symmetric spaces as well known. In all cases except for one, we will
show that they indeed carry a metric with \pc\ as a consequence of
\pref{fib}.
\smallskip

1) The first class of examples are the homogeneous flag manifolds
due to Wallach: $W^6=\SU(3)/\T^2$, $W^{12}=\Sp(3)/\Sp(1)^3$ and
$W^{24}=\F_4/\Spin(8)$. They are the total space of the following
homogeneous fibrations:
\bigskip

 \centerline{$\Sph^2\to \SU(3)/\T^2 \to \CP^2$,}

\bigskip

 \centerline{$\Sph^4\to \Sp(3)/\Sp(1)^3 \to \HP^2$,}

\bigskip

 \centerline{$\Sph^8\to \F_4/\Spin(8) \to \CaP$.}

\bigskip

\no  We now show that $W^6=\SU(3)/\T^2$ has \pc, the other cases
being similar.  Consider the inclusions
$\T^2\subset\U(2)\subset\SU(3)$ giving rise to the above homogeneous
fibration. Here we embed $\U(2)$ as the upper $2\times 2$ block,
i.e., $\U(2)=\{\diag(A,\det \bar{A})\mid A\in\U(2)\}$. A vector in
$\fm$ is of the form $Y=\left(
                                   \begin{array}{cc}
                                     0 & v  \\
                                     -\bar{v} &  0 \\
                                   \end{array}
                                 \right)$
with $v\in\C^2$ and one easily shows that $[A,Y]=Av+\tr(A)v$ for
$A\in\fu(2)$. Hence if $X\in\fp\subset\fsu(2)$ and $Y\in\fm$,
$[X,Y]=0$ iff $X=0$ or $Y=0$. This shows that part (c) of \pref{fib}
holds. As for $(a)$ and $(b)$ the  fiber and base are symmetric
spaces of rank 1 and thus $\SU(3)/\T^2$ has \pc. On the other hand,
one easily sees that there are vectors $X,Y\in\fm+\fp$ with
$[X,Y]=0$. Thus the biinvariant metric has \nnc\ but with some
0-curvature 2-planes. The Cheeger deformation deforms this metric
into one with \pc.

\smallskip

2) The Berger space $B^{13}=\SU(5)/\Sp(2)\cdot \S^1$  admits a
fibration

\bigskip
 \centerline{$\RP^5\to\SU(5)/\Sp(2)\cdot \S^1\to\CP^4,$}
\bigskip

\no coming from the inclusions $\Sp(2)\cdot
\S^1\subset\U(4)\subset\SU(5)$. Here $\Sp(2)\subset\SU(4)$ is the
usual embedding and $\S^1$ is the center of $\U(4)$. Furthermore,
the fiber is $\U(4)/\Sp(2)\cdot \S^1=\SU(4)/\Sp(2)\cdot \Z_2
=\SO(6)/\O(5)=\RP^5$. A biinvariant metric on $\SU(5)$ restricts to
a biinvariant metric on $\SO(6)$ which induces a metric  with
constant curvature on the fiber $\RP^5$. The base is clearly a
symmetric space of rank 1 and condition (c) is  verified as in the
previous case.

\smallskip

3) The Aloff-Wallach spaces
$W^7_{p,q}=\SU(3)/\diag(z^{p},z^{q},\bar{z}^{p+q}) , (p,q)=1,$ form
an infinite family.  We claim that they have positive curvature iff
$pq(p+q)\ne 0$. They admit a fibration

\bigskip
\centerline{$\Sph^3/\Z_{p+q}\to W_{p,q}\to \SU(3)/\T^2,$}
\bigskip

\no coming from the inclusions
$\diag(z^{p},z^{q},\bar{z}^{p+q})\subset\U(2)\subset\SU(3)$. Hence,
as long as $p+q\ne 0$, the fiber is the lens space
$\U(2)/\diag(z^{p},z^{q})=\SU(2)/\diag(z^{p},z^{q})$ with
$z^{p+q}=1$.

A vector in $\fm$ again  has the form $Y=\left(
                                   \begin{array}{cc}
                                     0 & v  \\
                                     -\bar{v} &  0 \\
                                   \end{array}
                                 \right)$.
Since the Lie algebra of $H$ is spanned by $\diag(ip,iq,-(ip+iq))$,
we can write an element in $\fp$ as $X=\diag(A,-\tr A)$ where
$A=\left(
     \begin{array}{cc}
       i(2q+p)a & z \\
       -\bar{z} & -i(q+2p)a \\
     \end{array}
   \right) $
   with $a\in\R$ and $z\in\C$. Hence $[X,Y]=Av+i(q-p)av$, i.e., $A$
   has an eigenvalue $i(p-q)a$ if $v\ne 0$. But one easily shows that
    this is only possible when $A=0$ or $pq=
   0$. Hence condition (c)
   is satisfied and since $(a)$ and $(b)$ clearly hold, $W_{p,q}$
   has positive curvature as long as $pq(p+q)\ne 0$. If on the other
   hand one of $p,q$ or $p+q$ is 0, say $p=0$,  one easily shows that the fixed point
   set of $\diag(1,-1,-1)\in H=\diag(z^{p},z^{q},\bar{z}^{p+q})=\diag(1,z,\bar{z}) $ is equal to
   $\U(2)/\diag(z,\bar{z})=\Sph^2\times\Sph^1/\Z_2$. Since fixed point sets of isometries
   are totally geodesic, and since $\Sph^2\times\Sph^1/\Z_2$ does
   not carry a metric with \pc, these Aloff-Wallach spaces cannot carry a
   homogeneous metric with \pc.

\smallskip

\smallskip

4. Finally  we have the Berger space: $B^7=\SO(5)/\SO(3)$. To
describe the embedding $\SO(3)\subset\SO(5)$, we recall that
$\SO(3)$ acts orthogonally via conjugation on the set of $3\times 3$
symmetric traceless matrices. This space is special since $\SO(3)$
is maximal in $\SO(5)$ and hence does not admit a homogeneous
fibration. It is also what is called isotropy irreducible, i.e., the
isotropy action of $H$ on the tangent space is irreducible. This
implies that there is only one $\SO(5)$ invariant metric up to
scaling. Now a direct calculation is necessary in order to show that
a biinvariant metric on $\SO(5)$ induces \pc\ on $B^7$.

\smallskip

\begin{rem*}
a) The examples $B^7$ and $B^{13}$ were found by Berger
\cite[\!\!'61]{Be1} when classifying  normal homogeneous metrics
with \pc\, i.e., metrics on $G/H$ induced by a biinvariant metric on
$G$.  But in \cite['99]{Wi1}) B.Wilking observed that the
Aloff-Wallach space $W_{1,1}$ is missing since it can  be written as
$\SU(3)\SO(3)/\U(2)$ where a biinvariant metric induces positive
curvature.

\smallskip

b)In \cite[\!\!'72]{Wa} Wallach also  proved that if one adds to the
assumptions in \pref{fib} that the fiber is a  symmetric pair as
well, then  the metrics $g_t$ with $1<t<4/3$ have positive curvature
also. This applies to the flag manifolds and the Berger space
$B^{13}$. We do not know of a simple geometric proof of this fact,
similar to the one we gave in \pref{fib}. It is also mysterious that
the limiting value $4/3$ is the same as in \lref{leftinv}. The
number $4/3$ shows up again if one considers homogeneous metrics on
spheres, scaled in the direction of one of the Hopf fibrations with
fibers $\Sph^1, \Sph^3$ or $\Sph^7$. As was shown in \cite['07]{VZ},
they have positive curvature as long as the scale is less than
$4/3$. In the cases where the fiber is 3 or 7 dimensional, the proof
again requires detailed curvature estimates. It would be interesting
to obtain a uniform and less computational understanding why the
number 4/3 appears in all 3 cases.
\end{rem*}

After a classification of all homogeneous spaces  which admit a
metric with positive curvature, one can ask for the {\it best}
homogeneous metric, i.e., the one with largest pinching $\delta$.
This is a rather difficult question since pinching constants are
notoriously difficult to compute. For the homogeneous spaces which
are not symmetric spaces of rank 1 this was done in
\cite[\!\!'79]{Va} for the flag manifolds and in \cite[\!\!'99]{Pu}
for the remaining cases (see also \cite[\!\!'66]{El},
\cite[\!\!'71]{He}, \cite[\!\!'81]{Hu} for previous work).
Interestingly, one obtains three homogeneous spaces, $B^7,B^{13}$
and $W_{1,1}$ which admit metrics with pinching $\delta=1/37$. In
the first two cases this is the best metric and in the latter case
the best one among all metrics invariant under $\SU(3)\SO(3)$. For
the flag manifolds the best metric has pinching $\delta=1/64$. In
\cite[\!\!'99]{Pu} one finds numerical values for the pinching
constants of the best homogeneous metrics on $W_{p,q} \ne W_{1,1}$.
It turns out to be an increasing function of $p/q$ when $0<p\le q$
and is in particular always $< 1/37$. For $W_{1,1}$ the set of
$\SU(3)$ invariant metrics is 10-dimensional, which makes the
computations rather difficult, even numerically.

\bigskip

\centerline{\it Biquotients with \pc.}

\bigskip

As explained in Section 2, biquotients $G/\!/H$ are obtained when
$H\subset G\times G$ acts on $G$ from the left and from the right.
When the action is free, the biinvariant metric on $G$ induces a
metric on $G/\!/H$ with \nnsc. In some cases, this can be deformed
via a Cheeger deformation into one with \pc.   We now describe these
biquotient examples explicitly and prove that they have positive
curvature.

\smallskip

1) There is an analogue of the 6-dimensional flag manifold which is
a biquotient of $\SU(3)$ under an action of $\T^2=\{(z,w)\mid
z,w\in\C \ , \ |z|=|w|=1\}$. It is given by:
$$E^6=\diag(z,w,zw)\backslash\SU(3)/\diag(1,1,z^2w^2)^{-1}.$$
The action by $\T^2$ is clearly free. In order to show that  this
manifold is not diffeomorphic  to the homogeneous flag $W^6$, one
needs to compute the cohomology with integer coefficients. The
cohomology groups are the same for both manifolds, but the ring
structure is different (\cite[\!\!'84]{E2}). The fact that this
manifold admits a metric with \pc\ will follow from the next
example.

\bigskip

2) We now describe the 7-dimensional family of Eschenburg spaces
$E_{k,l}$, which can be considered as a generalization of the Aloff
Wallach spaces. Let $k:= (k_1,k_2,k_3)$ and $l:= (l_1,l_2,l_3) \in
\Z^3$ be two triples of integers with $\sum k_i=\sum l_i $. We can
then define a two-sided action of $\S^1 = \{z\in \C \mid |z|=1 \}$
on $\SU(3)$ whose quotient we denote by $E_{k,l}$:
$$E_{k,l} :=  \diag(z^{k_1}, z^{k_2},
z^{k_3})\backslash \SU(3)/ \diag(z^{l_1}, z^{l_2}, z^{l_3})^{-1}.$$
\no The action is free if and only if
$\diag(z^{k_1},z^{k_2},z^{k_3})$ is not conjugate to
$\diag(z^{l_1},z^{l_2},z^{l_3})$, i.e.
\begin{equation*}\label{free}
\gcd(k_1-l_{i}\ , k_{2}-l_{j})=1,\ \ \text{\it for all }\ \  i\neq
j\,,\,i,j \in \{1,2,3\}\,.
\end{equation*}

We now claim:
\begin{prop}\label{Epos}
An Eschenburg space $E_{k,l}$ has \pc\ if
\begin{equation*}\label{pos}
k_i \notin [\min(l_1,l_2,l_3),\max(l_1,l_2,l_3)]
\end{equation*}
holds for all $1\leq i \leq 3$.
\end{prop}

\begin{proof}
As a  metric we choose the one  induced by a left invariant metric
on $\SU(3)$, in fact the same one as in \lref{sym} that we used for
$W^6$ and $W^7_{p,q}$. We first describe in a more explicit fashion
the set of 0-curvature 2 planes.

\begin{lem}\label{su3}
Let $Q_t$ be a left invariant metric on $\SU(3)$ as  in  \lref{sym}
with $G=\SU(3)$ and $K=\U(2)=\diag(A,\det \bar{A})$.  A 0-curvature
2-plane either contains a vector of the form $X=\diag(i,i,-2i)$,
which lies in the center of $\U(2)$, or one of the form
$X=\Ad(k)\diag(-2i,i,i)$ for some $k\in K$.
\end{lem}
\begin{proof}
By \lref{sym} a 0-curvature 2-plane is spanned by $X,Y$ with
$[X,Y]=[X_{\fk},Y_{\fk}]=[X_{\fp},Y_{\fp}]=0$. Since
$X_{\fp},Y_{\fp}$ are tangent to $G/K=\CP^2$, they   are linearly
dependent, and we can thus assume that $X_{\fp}=0$. If $X,Y$ both
lie in $\fk$, the fact that $[X,Y]=0$ implies that the 2-plane
intersects the center of $\fu(2)\cong\R\oplus\fsu(2)$, i.e. it
contains $X=\diag(i,i,-2i)$. If not, let $X=\diag(A,-\tr A)$ and
$Y_{\fp}=\left(
                                   \begin{array}{cc}
                                     0 & v  \\
                                     -\bar{v} &  0 \\
                                   \end{array}
                                 \right)$ with $0\ne v\in\C^2$. Then
                                 $0=[X,Y_{\fp}]=Av+(\tr A)v$ implies
                                 that $-\tr A$ and $2\tr A$ are
                                 eigenvalues of $A$ which means $A$
                                 is conjugate to $\diag(-2i,i)$,
                                 which proves our claim.
\end{proof}
In order to show that $Q_t$ induces positive curvature on $E_{k,l}$,
we need to prove that a 0-curvature 2-plane can never be horizontal,
i.e., it cannot be orthogonal to the vertical direction of the
$\S^1$ action.      Let $X_1=i\diag(k_1,k_2,k_3)$ and
$X_2=i\diag(l_1,l_2,l_3)$. Then the vertical space at $g\in\SU(3)$
is spanned by $(R_g)_*(X_1)-(L_g)_*(X_2)$, where $R_g$ and $L_g$ are
right and left translations. Since the metric is left invariant, we
can translate horizontal and vertical space back to $e\in\SU(3)$ via
$L_{g^{-1}}^*$. Thus the translated vertical space is spanned by
$\Ad(g^{-1})X_1-X_2$. We now need to show that a vector as in
\lref{su3} can never be orthogonal to it.

To facilitate this computation, observe the following. If $\ft$ is
the Lie algebra of a maximal torus in $G$, then critical points of
the function $g\to Q(\Ad(g)A,H)$ for fixed $A,H\in\ft$ are obtained
when $\Ad(g)A\in\ft$ also. Indeed, if $g_0$ is critical, we have
$0=Q([Y,\Ad(g_0)A],H)=Q([\Ad(g_0)A,H],Y)$ for all $Y\in\fg$ and thus
$[\Ad(g_0)A,H]=0$. For a generic vector $H\in\ft$ we have that
$\exp(tH)$ is dense in the compact torus $\exp(\ft)$ and hence
$[\Ad(g_0)A,\ft]=0$ which by maximality of $\ft$ implies that
$\Ad(g_0)A\in\ft$. If $H$ is not generic, the claim follows by
continuity.

We now apply this to the function
$Q_t(\Ad(g)X_1-X_2,\diag(i,i,-2i))$ which we need to show is never
0. This amounts to showing that $Q(\Ad(g)X_1,\diag(i,i,-2i))\ne $
\linebreak $ Q(\diag(i,i,-2i) ,X_2)= l_1+l_2-2l_3$. But maximum and
minimum of the left hand side, according to the above observation,
lies among the values $k_r+k_s-2k_t$, $r,s,t$ distinct. Subtracting
$\sum k_i = \sum l_i$ we see that one needs to assume that
$l_3\notin [\min(k_i),\max(k_i)]$. Next, according to \lref{su3}, we
need $Q(\Ad(g)\diag(-2i,i,i),X_1) \ne Q(\Ad(k)\diag(-2i,i,i),X_2) )
$ for any $g\in G$ and $k\in K$. According to the above principle,
the left hand side has max and min among $k_r+k_s-2k_t$ whereas the
right hand side among $-2l_1+l_2+l_3 \ , \ -2l_2+l_1+l_3$. Thus we
need to assume that the interval $[\min(l_1,l_2),\max(l_1,l_2)]$
does not intersect $[\min(k_i),\max(k_i)]$. This is one of the
possible cases. To obtain one of the other ones, we can choose a
different block embedding for $K=\U(2)\subset\SU(3)$.
\end{proof}

\smallskip

Among the biquotients $E_{k,l}$ there are two interesting
subfamilies. $E_p=E_{k,l}$ with $k=(1,1,p)$ and $l=(1,1,p+2)$ has
positive curvature when $p>0$. It admits a large group acting by
isometries. Indeed, $G=\SU(2)\times\SU(2)$ acting on $\SU(3)$ on the
left and on the right, acts by isometries in the Eschenburg metric
and commutes with the $\S^1$ action. Thus it acts by isometries on
$E_p$ and one easily sees that $E_p/G$ is one dimensional, i.e.,
$E_p$ is \coo. A second family consists of the \co\ two Eschenburg
spaces $E_{a,b,c}=E_{k,l}$ with $k=(a,b,c)$ and $l=(1,1,a+b+c)$.
Here $c=-(a+b)$ is the subfamily of Aloff-Wallach spaces. The action
is free iff $a,b,c$ are pairwise relatively prime and the Eschenburg
metric has \pc\ iff, up to permutations, $a\ge b\ge c >0$ or $a\ge b
>0\ , c<-a$.
For these spaces $G=\U(2)$ acts by isometries on the right and
$E_{a,b,c}/G$ is two dimensional. For a general Eschenburg space
$G=\T^3$ acts by isometries and $E_{k,l}/G$ is four dimensional. In
\cite[\!\!'06]{GSZ} it was shown that these groups $G$ are indeed
the id component of the full isometry group of a positively curved
Eschenburg space.

\smallskip

To see that the biquotient $\SU(3)/\!/T^2$ has positive curvature,
we can view it as an $\S^1$ quotient of the Eschenburg spaces
$\diag(z^p,z^q,z^{p+q})\backslash \SU(3)/\diag(1,1,z^{2p+2q})^{-1}$
which has positive curvature when $pq>0$.

\bigskip

3) We finally have the 13-dimensional Bazaikin spaces $B_{q}$, which
can be considered as a generalization of the Berger space $B^{13}$.
Let $q=(q_1,\dots,q_5)$ be a 5-tuple of integers with ${q}=\sum q_i$
and define
\begin{equation*}\label{Ba2}
B_q = \diag (z^{q_1} , \dots , z^{q_5} ) \backslash \SU(5)/ \diag
(z^{{q}},A)^{-1},
\end{equation*}
where $A\in \Sp(2)\subset \SU(4)\subset\SU(5)$. Here we follow  the
treatment in \cite[\!\!'98]{Zi1} of Bazaikin's work
\cite[\!\!'96]{Ba1}. First, one easily shows that the action of
$\Sp(2)\cdot\S^1$ is free if and only if
\begin{equation*}\label{e:gcd}
\text{\it all $q_i$'s are odd\ \ and\ \ }\gcd(q_{\gs (1)}+q_{\gs
(2)},q_{\gs (3)}+ q_{\gs (4)})=2,
\end{equation*}
for all permutations $\gs\in S_5.$ On $\SU(5)$ we choose an
Eschenburg metric by scaling the biinvariant metric on $\SU(5)$  in
the direction of $\U(4)\subset\SU(5)$. The right action of
$\Sp(2)\cdot\S^1$ is then  by isometries. Repeating the same
arguments as in the previous case, one shows that the induced metric
on $\SU(5)/\!/\Sp(2)\cdot\S^1$ satisfies
\begin{equation*}\label{pos}
\sec>0 \text{\it \quad if and only if\quad }  q_i+q_j>0 \ (\text{or}
<0)\; \text{\it for all }  i < j .
\end{equation*}
The special case of $q=(1,1,1,1,1)$ is the homogeneous Berger space.
One again has a one parameter subfamily that is \coo, given by
$B_p=B_{(1,1,1,1,2p-1)}$ since $\U(4)$ acting on the left induces an
isometric action on the quotient. It has \pc\ when $p\ge 1$.

\bigskip

Unlike in the homogeneous case, there is no general classification
of positively curved biquotients, except  in the following cases. We
call a  metric on $G /\!/H$ torus invariant if it is induced by a
left invariant metric on $G$ which  is also right invariant under
the action of a maximal torus. The main theorem in
\cite[\!\!'84]{E2} states that an even dimensional biquotient $G
/\!/H$ with $G$ simple and  which admits a positively curved torus
invariant metric is diffeomorphic to a rank one symmetric space or
the biquotient $\SU(3)/\!/T^2$. In the odd dimensional case he shows
that $G /\!/H$ with a positively curved torus invariant metric and
$G$ of rank 2 is either diffeomorphic to a homogeneous space or a
positively curved Eschenburg space. In particular, the sufficient
conditions in \pref{Epos} are also necessary not only for Eschenburg
metrics, but more generally torus invariant metrics. The
classification of the remaining odd dimensional biquotients with
$\rank G >2$ was taken up again in \cite[\!\!'98]{Bo} where it was
shown that  if one assumes in addition that $H=H_1\cdot H_2$ with
$H_1=\S^1$ or $H_1=\S^3$  and such that $H_2$ contains no rank one
normal subgroups and operates only on one side of $G$, the manifold
is diffeomorphic to a homogeneous space, an Eschenburg spaces, or a
Bazaikin space with positive curvature. The case where $G$ is not
simple, on the other hand, is wide open. As we will see in Section
5, one obtains a large number of examples with almost positive
curvature in this more general class of biquotients.

\smallskip

Not much is known about the pinching constants of the positively
 curved  metrics on biquotients. One easily sees that for
 a sequences of Eschenburg spaces $E_{k,l}$ where $(k/|k|, l/|l|)$
 converges to $((1,1,-2)/\sqrt{6},(0,0,0))$,  the pinching  of the Eschenburg metric
converges
 to $1/37$. In \cite[\!\!'04]{Di} W.Dickinson  proved that for a general
 positively curved
 Eschenburg space $E_{k,l}$ with its Eschenburg metric one has $\delta\le
 1/37$ with equality only for $W_{1,1}$. Furthermore, the
 pinching constant for the \coo\ Eschenburg spaces $E_p$ goes to $0$ when
 $p\to\infty$.

\bigskip

\centerline{\it Fundamental groups of positively curved manifolds}

\bigskip

A classical conjecture of S.S.Chern states that, analogously to the
Preismann theorem for negative curvature, an abelian subgroup of the
fundamental group of a positively curved manifold is cyclic. This is
in fact not true. The first counter examples were given in
\cite[\!\!'98]{Sh}, and further ones in \cite[\!\!'00]{GS} (see also
\cite[\!\!'99]{Ba2}):

\begin{thm} The following groups act freely on a positively curved
manifold:
\begin{itemize}
\item[(a)] (Shankar) $\Z_2\oplus\Z_2$ on the Aloff Wallach space
$W_{1,1}$ and the \coo\ Eschenburg space $E_2$.
\item[(b)](Grove-Shankar) The group $\Z_3\oplus\Z_3$ on the Aloff Wallach space
$W_{p,q}$ if $3$ does not divide $pq(p+q)$.
\end{itemize}
\end{thm}

In the case of $W_{1,1}=G/H=\SU(3)/\diag(z,z,\bar{z}^2)$ this
follows since $N(H)/H=\U(2)/Z(\U(2))=\SO(3)$ acts isometrically in
the Eschenburg metric and the right action of $N(H)/H$ is  free on
any $G/H$. Thus a finite subgroup of $\SO(3)$ acts freely as well.
Further example  are given in  \cite[\!\!'06]{GSZ} for the
cohomogeneity two Eschenburg spaces. But there are no examples known
where $\pi_1(M)=\Z_p\oplus\Z_p$, $p>3$ a prime, acts freely on a
positively curved manifold.

\bigskip

\centerline{\it Topology of positively curved examples}

\bigskip

In dimension 7 and 13 we have infinitely many homotopy types of
positively curved manifolds since an Eschenburg spaces satisfies
$H^4(E_{k,l},\Z)=\Z_r$ with $r=\sigma_2(k)-\sigma_2(l)$
 and for a
Bazaikin spaces one has $H^6(B_q,\Z)=\Z_r$ with
$8r=\sigma_3(q)-\sigma_1(q)\sigma_2(q)$ where $\sigma_i$ is the
elementary symmetric polynomial of degree $i$. On the other hand,
for a fixed cohomology ring, there are only finitely many known
examples \cite[\!\!'06]{CEZ},\cite[\!\!'06]{FZ1}. A classification
of 7-dimensional manifolds whose cohomology type is like that of an
Eschenburg space was obtained by Kreck-Stolz \cite[\!\!'91]{KS1} in
terms of certain generalized Eells-Kuiper invariants. They also
computed these invariants for the Aloff Wallach spaces and obtained
examples that are homeomorphic but not diffeomorphic. Kruggel
\cite[\!\!'05]{Kr} computed the Kreck-Stolz invariants for a general
Eschenburg space in terms of number theoretic sums and
Chinburg-Escher-Ziller \cite[\!\!'06]{CEZ} found further examples of
this type.

\begin{thm} One has the following examples with positive curvature:
\begin{itemize}
\item[(a)]  The pair of Aloff Wallach spaces $W_{k,l}$ with
$(k,l)=(56.788\; , \; 51.561)$ and $(k,l)=(61.213\; ,\; 18.561 )$
and the pair of Eschenburg spaces $E_{k,l}$ with $(k\, ;\, l)=(79,
\,49, \,-50 \ ; $ $ 0, \,46, \,32)$ and $(k \, ;\, l)=(75, \,54,
\,-51 \ ;\, 0, \,46 , \,32 )$ are homeomorphic to each other but not
diffeomorphic.
\item[(b)] The
pair of Aloff Wallach spaces $W_{k,l}$ with $(k,l)=(4.638.661\; ,\;
4.056.005 )$ and $(k,l)=( 5.052.965\; ,\; 2.458.816 )$ and  the pair
of Eschenburg spaces $E_{k,l}$ with $(k\, ;\, l)=(2.279\, ,
\,1.603\, , \,384 \ ;\, 0,0, \, 4.266 )$ and $(k\, ;\, l)= (2.528\,
, \,939\, , \,799 \ ;\, 0 ,0 ,\, 4.266 )$ are diffeomorphic to each
other but not isometric.
\end{itemize}
\end{thm}

The diffeomorphic pairs of Aloff-Wallach spaces give rise to
different components of the moduli space of positively curved
metrics \cite[\!\!'93]{KS2}, in fact these two metrics cannot be
connected even by a path of metrics with positive scalar curvature.
The diffeomorphic pair of Eschenburg spaces are interesting since
such cohomogeneity two manifolds also carry a 3-Sasakian metric (see
 Section 6) and they give rise to the first known manifold which
carries two non-isometric 3-Sasakian metrics.

\smallskip

The situation for the Bazaikin spaces seems much more rigid. A
computation of the Pontryagin classes and the linking form seems to
indicate, verified for the first 2 Billion examples, that Bazaikin
spaces are all pairwise diffeomorphically distinct, see
\cite[\!\!'06]{FZ1}. There is also only one Bazaikin space, the
Berger space $B^{13}$, which is homotopy equivalent to a homogeneous
space.

\smallskip

The Berger space $B^7=\SO(5)/\SO(3)$ plays  a special role. It is,
apart from spheres, the only known odd dimensional positively curved
manifold which is 2-connected, which should be compared with the
finiteness theorem by Fang-Rong and Petrunin-Tuschmann mentioned in
Section 1. It is also, apart from the Hopf bundle, the only $\Sph^3$
bundle over $\Sph^4$ which is known to have positive curvature since
it was shown in \cite[\!\!'04]{GKS} that it is diffeomorphic to such
a bundle. The topology of $\Sph^3$ bundles over $\CP^2$ is studied
in  \cite{EZ} where it is shown that they are frequently
diffeomorphic to positively curved Eschenburg spaces when the bundle
is not spin. It is thus natural to ask:

\begin{problem}
Does every $\Sph^3$ bundle over $\Sph^4$, and every $\Sph^3$ bundle
over $\CP^2$ which is not spin, admit a metric with \pc. Do
$\S^3$-principle bundles over $\Sph^4$, and $\SO(3)$-principle
bundles over $\CP^2$ which are not spin, admit a metric with \pc\
invariant under the principal bundle action.
\end{problem}

Notice that the existence in the latter case would imply, by
considering the associated $\Sph^2$ bundles, the existence of
infinitely many homotopy types of positively curved manifolds in
dimension 6. Also recall that there are two $\SO(3)$-principle
bundles over $\CP^2$ with such positively curved metrics and that
all bundles in Problem 6 have a metric with \nnc.

\bigskip

\section{Examples with almost positive or almost non-negative curvature}

\bigskip

As was suggested by Fred Wilhelm, there are two natural classes of
metrics that lie in between \nnc\ and \pc. In an initial step of
deforming a \nn ly curved metric into one with \pc\ one can first
make the curvature of all two planes at a point positive. We  say
that a metric has {\it \qp\ } curvature if there exists an open set
such that all \sc s  in this open set are positive. In a second step
one can try to deform the metric so that all \sc s in an open and
dense set are positive. We say that a metric with this property has
{\it \ap } \cu. It is natural to suggest that there should be
obstructions to go from \nn\ to \qp\ \cu\ and from \ap\ to \pc, but
that one should always be able to deform a metric from \qp\ to \ap\
\cu.

\smallskip

The first example of a manifold with almost positive curvature was
given by P.Petersen and F.Wilhelm in \cite[\!\!'99]{PW}, where it
was shown that $T_1\Sph^4$ has this property. In \cite[\!\!'01]{W}
it was shown that the Gromoll Meyer sphere
 $\Sigma^7=\Sp(2)/\!/\Sp(1)$ admits a metric with almost positive
 curvature as well. Following suggestions in \cite['84]{E2}, a
 simpler proof for a slightly different metric on $\Sigma^7$
  was given in \cite[\!\!'02]{E3} with a corrected proof in
  \cite[\!\!'07]{EK}.

\bigskip

We now describe some remarkable examples of metrics with \ap\ \cu\
due to Wilking \cite[\!\!'02]{Wi2}:

\begin{thm}[Wilking]\label{almost pos}
Let $M$ be one of the following manifolds:
\begin{itemize}
\item[(a)] One of the  projectivised tangent bundles $P_\R T(\RP^n)$,
$P_\C T(\CP^n)$, or $P_\QH T(\HP^n)$.

\item[(b)] The homogeneous space $M^{4n-1}_{p,q}=U(n+1)/H_{k,l}$ with
$H_{k,l}=\{\diag(z^p,z^q,A)\mid |z|=1, A\in\U(n-1)\}$, where $pq<0$
and $(p,q)=1$.
\end{itemize}
Then $M$ carries a metric with \ap\ \cu.
\end{thm}
Here projectivised means that we identify a tangent vector $v$ with
$\lambda v$ where $\lambda\in\R,\, \C$ or $\QH$ respectively. Notice
that in contrast to the known \pcu\ examples, these \ap ly curved
manifolds exist in arbitrarily high dimensions. Furthermore, in the
case of $n=2$ these are one of the known manifolds with positive
curvature, namely $P_\C T(\CP^2)$ and $P_\QH T(\HP^2)$ are the flag
manifolds $W^6$ and $W^{12}$ and $M^{7}_{p,q}$ is the Aloff Wallach
space $W_{p,q}$. Notice that the unique Aloff Wallach space
$W_{1,0}=W_{1,-1}$ which does not admit a homogeneous metric with
\pc, thus admits a metric with \ap\ \cu.

These examples also show that in general a metric with almost
positive curvature  cannot be deformed to positive curvature
everywhere. Indeed, $P_\R T(\RP^{2n+1})$ is an odd dimensional
non-orientable manifold and hence by Synge's theorem does not admit
\pc. A particularly interesting special case is $P_\R
T(\RP^{3})=\RP^3\times\RP^2$ and $P_\R T(\RP^{7})=\RP^7\times\RP^6$.
If the manifold is compact and simply connected, it is not known
whether an almost positively curved metric can be deformed to
positive curvature. On the other hand,   there either are
obstructions or the generalized Hopf conjecture on
$\Sph^3\times\Sph^2$ is false.

\begin{proof}
We prove \tref{almost pos} in the simplest case of $P_\R
T(\RP^3)=\RP^3\times\RP^2$.

Define a left invariant metric $g$ on $G=\S^3\times\S^3$ by scaling
a biinvariant metric on $\fg$ in the direction of the diagonal
subgroup $K=\Delta \S^3$ as in \lref{sym}. Since $G/K=\S^3$ is a
symmetric pair of rank one, and since $(X,Y)_\fk=\frac 1 2
(X+Y,X+Y)$ , $(X,Y)_\fm=\frac 1 2 (X-Y,-X+Y)$, \lref{sym} implies
that a 0 curvature plane is spanned by vectors $(u,0) , (0,u)$ with
$0\ne u \in \Im \QH$. Here we regard $\S^3$ as the unit quaternions
with Lie algebra $\Im \QH$. $G$ acts on $T_1\Sph^3=\{(p,v)\mid
|p|=|v|=1 , \langle u,v\rangle =0\}$ via $(q_1,q_2)\star (p,v)=(q_1
p q_2^{-1}, q_1 v q_2^{-1})$ and the isotropy group of $(1,i)$ is
equal to $H=(e^{i\theta},e^{i\theta})$ and thus $G/H=T_1\Sph^3$. We
can rewrite the homogeneous space $G/H$ as a biquotient $\Delta G
\backslash G\times G/(1\times H)$ since $\Delta G \backslash G\times
G=G$.

We now claim that the product metric $g + g $ on $G\times G$ induces
a metric with almost positive curvature on $\Delta G \backslash
G\times G/(1\times H)$. For this, notice that each orbit of $\Delta
G$ acting on the left on $G\times G$ contains points of the form
$p=(\bar{a},\bar{b},1,1)$, $a,b\in\S^3$. The vertical space,
translated to the identity via left translation with $(a,b,1,1)$, is
equal to the direct sum of $(Ad(a)v,Ad(a)w,v,w)$ with $
v,w\in\Im\QH$ and $(0,0,i,i)\cdot \R$. If we set $g(A,B)=Q(PA,B)$
where $Q$ is a biinvariant metric on $G$, a horizontal vector is of
the form $\left( P^{-1}(-Ad(a)v,-Ad(b)w) , P^{-1}(v,w)\right)$ with
$(v,w)\bot (i,i)$. Since $P$ clearly preserves 2-planes spanned by
$(v,0),(0,v)$, a horizontal 0-curvature plane is spanned by $\left(
P^{-1}(-Ad(a)v,0) ,\right. $ $\left. P^{-1}(v,0)\right)$ and $\left(
P^{-1}(0,-Ad(b)v) , P^{-1}(0,v)\right)$ with   $Ad(a)v=\pm Ad(b)v$.
Thus $\bar{a}b$ either commutes or anticommutes with $v\in\Im\QH$
and since also $v\bot i$, either $\bar{a}b\bot i$ or $\bar{a}b\bot
1$. Hence points with 0-curvature 2 planes lie in two hypersurfaces
in $G/H$.

Since the group $L$ generated by $ (1,-1) $ and $(j,j)$ normalizes
$H$, and since the left invariant metric $g$ is also right invariant
under $L$, the quotient  $G/H\cdot L$ inherits a metric with almost
positive curvature and one easily sees that this quotient is $P_\R
T(\RP^3)=\RP^3\times\RP^2$.
\end{proof}

See \cite[\!\!'02]{Wi2} for further examples. In \cite[\!\!'03]{Ta2}
K.Tapp proved that the unit tangent bundles of $\CP^n$, $\HP^n$ and
$\CaP$, as well as the manifolds in \tref{almost pos} (b) with $pq >
0$, have quasi positive curvature. In \cite{Ke} M.Kerin showed that
all Eschenburg spaces have a metric with quasi positive curvature
and that $E_0$, the unique \coo\ Eschenburg space which does not
admit a \co\ metric with positive \cu, admits a metric with almost
positive curvature.

\smallskip

All known examples of  \ap\ \cu\ (and in fact all positively curved
examples as well) can be described, after possibly enlarging the
group, as so called {\it normal biquotients}, i.e., $M=G/\!/H$ with
metric on $M$ induced by a biinvariant metric on $G$. B.Wilking
showed in \cite[\!\!'02]{Wi2} that for such normal biquotients the
exponential image of a 0-curvature 2-plane  is totally geodesic and
flat. As was observed by K.Tapp \cite[\!\!'07]{Ta3}, this remains
true more generally for a Riemannian submersion $G\to M$ when $G$ is
equipped with a biinvariant metric. It is a natural question to ask
if the existence of an immersed flat 2-torus is sometimes an
obstruction to deform a metric with \nnc\
 to one with \pc.

 \smallskip

 In  \cite[\!\!'02]{Wi2} one finds a number of natural open
 questions:

\begin{itemize} \item  Can every quasi positively curved metric be
deformed to almost positive curvature.
\item Can a quasi positively curved metric where
the points with 0-curvatures are contained in a contractible set be
deformed to positive curvature.
\item Does an even dimensional almost positively curved manifold
have positive Euler characteristic.
\end{itemize}

\bigskip

\centerline{\it Almost non-negative curvature}

\bigskip

We say that a manifold $M$ has almost non-negative  curvature if
there exists a sequence of metrics $g_i$, normalized so that the
diameter is at most 1, with $\sec(g_i)\ge - 1/i$ for all $i\in\N$.
This includes the almost flat manifolds where $\sec(g_i)\le 1/i$ as
well. By  Gromov's almost flat manifold theorem, the latter are
finitely covered by a compact quotient of a nilpotent Lie group
under a discrete subgroup.

This is a much larger class of manifolds. Besides being invariant
under taking products, it is also well behaved under fibrations. In
\cite[\!\!'92]{FY}, Fukaya-Yamaguchi showed that:

\begin{thm}[Fukaya-Yamaguchi]
The total space of a principal $G$-bundle with $G$ compact  over an
almost non-negatively curved manifold is  almost non-negatively
curved as well.
\end{thm}

To see this, one puts a metric on the total space $M$ such that the
projection onto the base is a Riemannian submersion with totally
geodesic fibers. Scaling the metric on $M$ in the direction of the
fibers then has the desired properties as the scale goes towards 0.

Thus every associated bundle  $P\times_G F$, where $G$ acts
isometrically on a non-negatively curved manifold F,  has almost
non-negative curvature as well. This applies in particular to all
sphere bundles.

 As was shown in \cite[\!\!'04]{ST}, this class also
includes all \com s:

\begin{thm}[Schwachh\"ofer-Tuschmann]
Every compact \com\ has an invariant metric with  almost
non-negative sectional curvature.
\end{thm}
\begin{proof}
As was observed  by B.Wilking, this follows easily by using a
Cheeger deformation. If a compact group $G$ acts by \coo\ on $M$
then
 $M/G$ is either a circle or an interval. In the first case
  $M$ carries a $G$-invariant metric with \nnc. In
the second case we choose a $G$-invariant metric $g$ on $M$ which
has \nnc\ near the two singular orbits. This is clearly possible
since a neighborhood of a singular orbit is a homogeneous disk
bundle $G\times_K D$. We now claim that a Cheeger deformation $g_t$
with respect to the action of $G$ on $M$ has \ap\ curvature as
$t\to\infty$. Setting $s=1/t$, the metric on $M=M\times G/\Delta G$
is induced by a metric of the form $g+ s Q$ on $M\times G$. Thus its
diameter is clearly bounded as $s\to 0$. On the regular part, we can
assume that our 2-plane is spanned by vectors $X+\alpha T$ and $Y$,
where
 $X,Y$ are tangent to a principal orbit and $T$ is a unit vector
 orthogonal to all principal orbits. As we saw in Section 2,
the horizontal lift of these vectors to $M\times G$ under the
Riemannian submersion $M\times G\to M$ are of the form
$(sP^{-1}(sP^{-1}+\Id)^{-1}X +\alpha T,-(sP^{-1}+\Id)^{-1}X)$ and
$(sP^{-1}(sP^{-1}+\Id)^{-1}Y ,-(sP^{-1}+\Id)^{-1}Y)$. Only the first
component can contribute to a negative curvature. The curvature
tensor of this component goes to 0 with $s$, of order 2 if
$\alpha\ne 0$ and of order 4 if $\alpha=0$. On the other hand,
$|(X+\alpha T)\wedge Y|^2_{g_t}$ goes to 0 with order 1 if
$\alpha\ne 0$ and order 2 if $\alpha=0$. Hence the negative part
goes to 0 as $s\to 0$.
\end{proof}

The main obstruction theorems for almost non-negative curvature are:

\begin{itemize} \item  (Gromov) The Betti numbers are
universally bounded in terms of the dimension.
\item
(Fukaya-Yamaguchi \cite[\!\!'92]{FY}) The fundamental group contains
a nilpotent subgroup of finite index.
\item
(Kapovitch-Petrunin-Tuschmann) \cite[\!\!'06]{KPT2}) A finite cover
is a nilpotent space, i.e. the action of the fundamental group on
its higher homotopy groups is nilpotent,
\end{itemize}

Notice that a compact quotient of a nilpotent non-abelian Lie group
has almost non-negative curvature, but does not admit a metric with
\nnc. On the other hand, for simply connected manifolds there are no
known obstructions which  could distinguish between almost \nn\ and
\nn\ (or even positive) curvature. As was suggested by K.Grove, it
is also natural to formulate the Bott conjecture more generally for
almost non-negatively curved manifolds.

\bigskip

\section{Where to look for new examples?}

\bigskip

There are two natural suggestions where one might look for new
examples with  \psc. The first is given by a structure that almost
all known examples share. They are the total space of a Riemannian
submersion. If one considers the more general class where the base
space of the submersion is allowed to be an orbifold, then all known
examples share this property, see \cite[\!\!'06]{FZ2}.

\bigskip

\centerline{\it Fiber bundles}

\bigskip

A.Weinstein \cite[\!\!'80]{We} considered fiber bundles $\pi\colon
M\to B$, where $\pi$ is a Riemannian submersions with totally
geodesic fibers. He called such a bundle  {\it fat} if  all
vertizontal curvatures, i.e. the curvature of a 2-plane spanned by a
horizontal and a vertical vector, are positive. This seemingly week
assumption already  places strong restrictions on the bundle. In
fact, one has (\cite[\!\!'81]{DR}:

\begin{thm} (Derdzinski-Rigas) Every $\Sph^3$ bundle over $\Sph^4$
which is fat is a Hopf bundle.
\end{thm}

This negative result seems to have discouraged the study of fat
fiber bundles until recently. On the other hand, as we saw in
\pref{fib},  most homogeneous examples of \pc\ are the total space
of a fat bundle. See \cite[\!\!'00]{Zi2} for a survey of what was
known about fatness up to that point.  In \tref{3sasak} we will see
that there are infinitely many $\Sph^3$ orbifold bundles over
$\Sph^4$ which are fat. It is thus natural and important to study
fat bundles in this more general category.

\smallskip

A natural class of  metrics  is given by a connection metric on a
principle $G$-bundle $\pi\colon P\to B$, sometimes also called a
Kaluza Klein metric. Here one chooses a principal connection
$\theta$, a metric $g$ on the base $B$, and a fixed biinvariant
metric $Q$ on $G$ and defines:
$$g_t(X,Y)=tQ(\theta(X),\theta(Y))+g(\pi_*(X),\pi_*(Y)). $$

The projection $\pi$ is then a Riemannian submersion with totally
geodesic fibers isometric to $(G,tQ)$.  Weinstein observed that the
fatness condition is equivalent to requiring  that the curvature
$\Omega$ of $\theta$ has the property that $\Omega_u=Q(\Omega,u)$ is
a symplectic 2-form on the horizontal space for every $u\in\fg$. If
$G=\S^1$, this is equivalent to the base being symplectic. If one
wants to achieve positive curvature on the total space, we need to
assume, in addition to the base having positive curvature, that
$G=\S^1$, $\SU(2)$ or $\SO(3)$. In \cite[\!\!'92]{CDR} a necessary
and sufficient condition for positive curvature of such metrics  was
given. The proof carries over immediately to the category of
orbifold principle bundles, which includes the case where the $G$
action on $P$ has only finite isotropy groups.

    \begin{thm}         \label{chavesderigas}
    (Chaves-Derdzinski-Rigas)
    A connection metric $g_t$ on aa  orbifold $G$-principle bundle with $\dim G\le 3$
   has positive curvature, for $t$
    sufficiently small, if and only if
    $$
    \left(\nabla_x\Omega_u\right)(x,y)^2 < |i_x\Omega_u|^2 k_B(x,y),
    $$

    \noindent for all linearly independent horizontal vectors
     $x,y$ and $0\ne u\in\fg$.
    \end{thm}

    Here $k_B(x,y)=g(R_B(x,y)y,x)$ is the unnormalized sectional
    curvature and $i_x\Omega_u\ne 0 $
    is precisely the above fatness condition. We call a principle connection with this property
    {\it hyperfat}.

    The simplest examples of hyperfat principal connections are
    given by the Aloff-Wallach  spaces  $W_{k,l}$,
    considered as a circle bundle over $\SU(3)/T^2$, since the fibers
    of a homogeneous fibration are  totally geodesic.
     As mentioned earlier,
    $W_{1,1}$ can also be considered as an $\SO(3)$ principle
    bundle over $\CP^2$
    which thus carries an $\SO(3)$ hyperfat connection
    (a fact first observed
     in \cite[\!\!'94]{Ch}).  In the orbifold
    category one has many more examples.
    Recall that a metric is called {\it 3-Sasakian} if $\SU(2)$
acts isometrically and almost freely with totally geodesic orbits of
curvature 1. Moreover, for $U$ tangent to the $\SU(2)$ orbits and
$X$ perpendicular to them,  $X\wedge U$ is required to be
 an eigenvector
of the curvature operator $\hat{R}$  with eigenvalue 1. In
particular the vertizontal curvatures are equal to 1, i.e., the
bundle is fat. This gives rise to a large new class of fat orbifold
principle bundles, see \cite[\!\!'99]{BG} for a survey. The
dimension of the base is a multiple of 4, and its induced (orbifold)
metric is quaternionic K\"{a}hler with positive scalar curvature.
One easily sees that the condition on the curvature operator is
equivalent to $\nabla_x\Omega_u= 0$ for all $u\in\fg$. Hence we
obtain the following Corollary, which was reproved later in
\cite[\!\!'04]{De1}:

\begin{cor}\label{PosSas}
A 3-Sasakian manifold has \psc, after the metric in the direction of
the $\SU(2)$ orbits is scaled down sufficiently, if and only if the
quaternionic K\"ahler quotient has \psc.
\end{cor}

In \cite[\!\!'66]{Be2} it was shown that a quaternionic K\"ahler
manifold of dimension $4n>4$ has \psc\ if and only if it is
isometric to $\HP^n$, which also holds for orbifolds. When the base
is 4-dimensional, quaternionic K\"ahler is equivalent to being self
dual Einstein and here an interesting new family of examples arises.
In \cite[\!\!'94]{BGM} it was shown that the \co\ two Eschenburg
spaces $E_{a,b,c}=\diag(z^a,z^b,z^c)\backslash
\SU(3)/\diag(1,1,z^{a+b+c})$, with $a,b,c$ positive and pairwise
relatively prime, carry a 3-Sasakian metric with respect to the
right action by $\SU(2)$. The quotients are weighted projective
spaces and Dearricott examined their sectional curvatures in
\cite[\!\!'05]{De2}:

\begin{cor} (Dearricott)\label{weighted}
The principal connection for the 3-Sasakian manifold $E_{a,b,c}$
with $0< a\le b \le c$ is hyperfat if and only if $c^2<ab$.
\end{cor}

Although the total space also carries an Eschenburg metric with \pc,
the projection to the base in that case does not have totally
geodesic fibers.

\smallskip

Since many of the known examples are also the total space of sphere
bundles, it is natural to study this category as well. A connection
metric on a sphere bundle   can be defined in terms of a metric
connection $\nabla$ on the corresponding vector bundle. It induces a
horizontal distribution on the sphere bundle and the fibers are
endowed with a metric of constant curvature. An analogue of
\ref{chavesderigas} for sphere bundles was proved in
\cite[\!\!'03]{Ta1}:

    \begin{thm}[Tapp] \label{kristhm}
A connection metric on an orbifold sphere bundle $E\to B$ has \pc,
for sufficiently small radius of the fibers, if and only if
$$ \langle(\nabla_xR^\nabla)(x,y)u,v\rangle^2 < |R^\nabla(u,v)x|^2
k_B(x,y), $$ \noindent for all linearly independent $x,y\in T_pB$
and $u,v\in E_p$, where we have set $\langle R(u,v)x,y\rangle
=\langle R(x,y)u,v\rangle$
    \end{thm}

    Notice that $R^\nabla(u,v)x\neq 0$ for all $u\wedge v\neq 0$, $x\neq 0$,
    means that  $(x,y)\to<R^\nabla(x,y)u,v>$ is nondegenerate
    for all $u\wedge v\neq 0$. But this is simply  Weinstein's fatness
    condition for the sphere bundle.

    Known examples are
   the homogeneous Wallach flag manifolds $W^6,W^{12},W^{24}$.
   The Aloff-Wallach spaces $W_{p,q}$ with $p+q=1$ are hyperfat $\Sph^3$
   bundles over $\CP^2$. Furthermore, the 3-Sasakian Eschenburg
   spaces $E_{a,b,c}$ in \cref{weighted}, where
   one of $a,b,c$ is even,
   are hyperfat $\Sph^3$ orbifold bundles and the associated $\Sph^2$ bundle
    $E_{a,b,c}/\S^1\to E_{a,b,c}/\SU(2)$ is hyperfat as well. If the
   base is a manifold, the condition $\nabla R=0$ is rather
   restrictive. For example, if the base is a symmetric space, it was shown
   in \cite[\!\!'01]{GSW} that the bundle must be homogeneous.
   Homogeneous fat fiber bundles were classified in
   \cite[\!\!'75]{BB1}. If the fiber dimension is larger than one,
    the base is symmetric. But if we assume in addition
   that the base has positive curvature, only the Wallach spaces
   remain.

   \bigskip

   The fiber bundle structure for most of the Eschenburg spaces and
   Bazaikin spaces do not have totally geodesic fibers. The Berger
   space $B^7=\SO(5)/\SO(3)$ is also the total space of an $\SU(2)$
   orbifold
   principal bundle over $\Sph^4$, but the fibers are again not
   totally geodesic. It is therefore also natural to examine warped
    connection metrics on the total space
   where the metric on the fiber is multiplied by a function on the
   base, see \cite[\!\!'03]{Ta1}, \cite[\!\!'05]{STT}. But notice that
   the known fibrations of the Eschenburg and Bazaikin spaces,
   with fiber dimension bigger than one, are not of this form either.
    Connection
   metrics with \nnc\ have also been studied in
   \cite[\!\!'90]{SW}, \cite[\!\!'95]{Ya} and
\cite[\!\!'05]{STT}.

\bigskip

\centerline{\it Cohomogeneity one manifolds}

\bigskip

A second natural class of manifolds where one can search for new
examples, especially in light of \tref{cohom}, are manifolds with
low cohomogeneity. Positively curved homogeneous spaces are
classified, so \com s are the natural next case to study.

There are many \coo\ actions on symmetric spaces of rank one. Among
the examples of positive curvature discussed in Section 4, one has a
number of other \coo\ actions. As mentioned there, the positively
curved Eschenburg spaces $$E^7_p = \diag(z,z,z^p) \backslash
\SU(3)/\diag(1,1,\bar{z}^{p+2}) , p\ge 1 ,$$ and the Bazaikin spaces
$$B^{13}_p = \diag(z,z,z,z,z^{2p-1})\backslash\SU(5)/\Sp(2)
\diag(1,1,1,1,\bar{z}^{2p+3}), p \ge 1,$$ admit \coo\ actions. Two
further examples are the Wallach space $W^7_{1,1} =
\SU(3)\SO(3)/\U(2)$ where $\SO(3)\SO(3)$ acts by \coo, and the
Berger space $B^7 = \SO(5)/\SO(3)$ with  $\SO(4)\subset\SO(5)$
acting by \coo.

In even dimensions, L.Verdiani \cite[\!\!'04]{Ve} classified all
positively curved \com s. Here only rank one symmetric spaces arise.
In odd dimensions, one finds a "preclassification" in
\cite[\!\!'06]{GWZ}. In dimension seven, a new family of candidates
arises: Two infinite families $P_k, Q_k$, $k\ge 1$, and one isolated
manifold $R$. The group diagram for $P_k$ is similar to the one
considered in \tref{principals4}:

$$
H=\Delta Q \subset  \{ (e^{i\gt},e^{i\gt})\cdot H \, ,
(e^{j(1-2k)\gt},e^{j(1+2k)\gt})\cdot H\}\subset \S^3\times \S^3,
$$
whereas the one for $Q_k$ is given by
$$
H=\{(\pm 1, \pm 1), (\pm i, \pm i)\}\subset \{
(e^{i\gt},e^{i\gt})\cdot H \, , (e^{jk\gt},e^{j(k+1)\gt})\cdot
H\}\subset \S^3\times \S^3.
$$

\no $R$ is similar to $Q_k$ with slopes $(3,1)$ on the left, and
$(1,2)$ on the right.

\begin{thm}(Verdiani, Grove-Wilking-Ziller)
A simply connected cohomogeneity one manifold $M$ with an invariant
metric of positive sectional curvature is equivariantly
diffeomorphic to one of the following:
\begin{itemize}
\item
An isometric action on a rank one symmetric space,
\item
     One of $E^7_p,  B^{13}_p$ or
$B^7$,
\item
One of the  $7$-manifolds $P_k, Q_k$, or $R$,
\end{itemize}
with one of the actions described above.
\end{thm}

The first in each sequence $P_k,Q_k$ admit an invariant metric with
positive curvature since $P_1=\Sph^7$ and $Q_1=W^7_{1,1}$.

 Among the \com s with codimension 2 singular orbits, which all
 admit \nnc\ by \tref{codim2},
are two families like the above $P_k$ and $Q_k$, but where the
slopes for $K_\pm$  are arbitrary. It is striking that in positive
curvature, with one exception, only the above slopes are allowed.
The exception is given by the positively curved \co one action on
$B^7$, where the isotropy groups are like those for $P_k$ with
slopes $(1,3)$ and $(3,1)$. In some tantalizing sense then, the
exceptional Berger manifold $B^7$ is associated with the $P_k$
family in an analogues way as the exceptional candidate $R$ is
associated with the $Q_k$ family. It is also surprising that all
non-linear actions in the above Theorem, apart from the Bazaikin
spaces $B_p^{13}$, are \co\ one under a group locally isomorphic to
$\S^3\times\S^3$.

\smallskip

These candidates also have interesting topological properties. $Q_k$
has the same cohomology groups as $E_k$ with $H^4(Q_k,\Z)=\Z_k$. The
manifolds $P_k$ are all 2-connected with $H^4(P_k,\Z)=\Z_{2k-1}$.
Thus it is natural to ask:

\begin{problem}
Are any of the manifolds $Q_k$, $k>1$, diffeomorphic to a positively
curved Eschenburg space? Are any of the manifolds $P_k$, $k>1$,
diffeomorphic to an $\Sph^3$ bundle over $\Sph^4$?
\end{problem}
Manifolds of this type are classified by their Kreck-Stolz
invariants, but these can be very difficult to compute in concrete
cases.

\smallskip

A somewhat surprising property that these candidates also have is
that they admit fat principal connections, in fact they admit
3-Sasakian metrics:

\begin{thm}\label{3sasak}
$P_k$ and $Q_k$ admit 3-Sasakian metrics which are orbifold
$\S^3$-principal bundles over $\Sph^4$ respectively $\CP^2$.
\end{thm}

This follows easily (\cite[\!\!'06]{GWZ}) from the celebrated
theorem due to Hitchin \cite[\!\!'93]{Hi}  that $\Sph^4$ admits a
family of self dual Einstein orbifold metrics invariant under the
\co\ one action by $\SO(3)$, one for each $k\ge 1$, which is smooth
everywhere except normal to one of the singular orbits where  it has
angle $2\pi/k$. One then shows that the induced 7-dimensional
3-Sasakian metric has no orbifold singularities, and by comparing
the isotropy groups of the \coo\ actions, it follows that they are
equivariantly diffeomorphic to $P_k$ and $Q_k$. Unfortunately, the
self dual Einstein metric on the base does not have \pc, unless
$k=1$, corresponding to the smooth metrics on $\Sph^4$ respectively
$\CP^2$. So \cref{PosSas} does not easily yield the desired metrics
of positive curvature on $P_k$ and $Q_k$.

\smallskip

Hence these candidates surprisingly have both features, they admit
\coo\ actions, and are also the total space of an orbifold principle
bundle. Both of these properties thus suggest concrete ways of
finding new metrics with \pc. We thus end with our final problem.

\begin{problem}
Do all manifolds $P_k\ , Q_k$and $R$ admit a \co\ metric with \pc?
\end{problem}

Notice that a positive answer for the manifolds $P_k$ would give
infinitely many homotopy types of positively curved 2-connected
manifolds. Thus the pinching constants $\delta_k$, for any metric on
$P_k$, would necessarily go to 0 as $k\to \infty$, and $P_k$ would
be the first examples of this type. It is natural to suggest that
the manifolds $E_p$, although not 2-connected, should have the same
property since the pinching constant for the Eschenburg metric goes
to $0$ as $p\to \infty$, .

\providecommand{\bysame}{\leavevmode\hbox
to3em{\hrulefill}\thinspace}

\end{document}